\documentclass[10pt]{article}
\usepackage[margin=2.5cm]{geometry}
\usepackage[utf8]{inputenc}
\usepackage{bussproofs}
\usepackage{mathrsfs}
\usepackage{amsmath}
\usepackage{amsthm}
\usepackage{amssymb}
\usepackage[dvipsnames]{xcolor}
\usepackage{stmaryrd}
\usepackage{mathtools}
\usepackage{calc}
\usepackage{cmll}
\usepackage{url}
\usepackage{amsmath}
\usepackage{longtable}
\usepackage{stmaryrd}
\usepackage{{xargs,ifthen,xstring}}
\usepackage{authblk}
\usepackage{txfonts}
\usepackage{array}
\newtheorem{theorem}{Theorem}

\newtheorem{definition}{Definition}
\newtheorem{proposition}{Proposition}

\newtheorem{corollary}{Corollary}
\newtheorem{conjecture}{Conjecture}

\title{Uniform validity of atomic Split rule\\in monotonic proof-theoretic semantics}

\author{Antonio Piccolomini d'Aragona\\

Eberhard Karls Universität Tübingen\\

\texttt{antonio.piccolomini-daragona@uni-tuebingen.de}}
\date{}

\begin{document}

\maketitle

\begin{abstract}
    Proof-theoretic semantics (PTS) is normally understood today as \emph{Base-Extension Semantics} (B-eS), i.e., as a theory of proof-theoretic consequence over atomic proof systems. Intuitionistic logic ($\texttt{IL}$) has been proved to be incomplete over a number of variants of B-eS, including a \emph{monotonic} one where \emph{introduction rules} play a prior role (miB-eS). In its original formulation by Prawitz, however, PTS consequence is not a primitive, but a \emph{derived} notion. The main concept is that of \emph{argument structure} valid relative to atomic systems and assignments of \emph{reductions} for eliminating generalised detours of inferences in non-introduction form. This is called \emph{Proof-Theoretic Validity} (P-tV), and it can be given in a monotonic and introduction-based form too (miP-tV). It is unclear whether, and under what conditions, the incompleteness results proved for $\texttt{IL}$ over miB-eS can be \emph{transferred} to miP-tV. As has been remarked, the main problem seems to be that the notion of argumental validity underlying the miB-eS notion of consequence is one where reductions are either forced to be \emph{non-uniform}, or \emph{non-constructive}. Building on some Prawitz-fashion incompleteness proofs for $\texttt{IL}$ based on the notion of (intuitionistic) construction, I provide in what follows a set of reductions which are surely uniform (however uniformity is defined) and constructive, and which make the atomic Split rule logically valid over miP-tV, thus implying the incompleteness of $\texttt{IL}$ over a Prawitzian (monotonic, introduction-based) framework strictly understood.
\end{abstract}


\section{Introduction}

In current investigations on completeness relative to \emph{Proof-Theoretic Semantics} (PTS), the latter is typically understood as what, after \cite{sandqvist}, is today called \emph{Base-Extension Semantics} (B-eS)---an overview of PTS as conceived here is in \cite{schroeder-heisterSE} while, for an overview of completeness on PTS, see \cite{piecha, schroederheisterrolf}. B-eS is a theory of consequence that qualifies as proof-theoretic because the ‘‘models" which consequence is evaluated on are (sets of) \emph{atomic proof systems}, i.e., (sets of) sets of rules which govern the deductive behaviour of the atomic part of one's language---see e.g. \cite{nascimentothesis, piechaschroeder-heisterbases, schroederheisterrolf, staffordsowhat}. The consequence relation itself is defined by standard induction, with two main alternatives available:

\begin{itemize}
    \item the relation may be either \emph{monotonic} or \emph{non-monotonic}, depending on whether the fact it holds on a given proof system is (resp. is not) preserved when extending the system, and
    \item the relation may be \emph{introduction-based} or \emph{elimination-based}, depending on whether the inductive clauses mirror the introduction (resp. elimination) rules of Gentzen's Natural Deduction calculi.
\end{itemize}
In what follows, I shall be concerned with (propositional) \emph{monotonic introduction-based} B-eS (miB-eS), with respect to which many completeness or incompleteness result have been proved. For what interests me here, especially important is an incompleteness theorem for the intuitionistic propositional calculus ($\texttt{IL}$), proved by de Campos Sanz, Piecha and Schroeder-Heister \cite{piechadecampossanz, piechadecampossanzschroederheisterpeirce, piechaschroeder-heisterdecampossanz, piechaschroeder-heister2019}.

In its original formulation by Prawitz \cite{prawitz1971, prawitz1973}, however, PTS cannot be really understood as a theory of consequence. Consequence is here only a \emph{derived} notion, whereas the prior concept is that of \emph{valid argument}. An argument is made up of two components: a \emph{proof-structure}, namely, a derivation-tree in Natural Deduction style where arbitrary inferences may occur, and a set of \emph{reductions} which, mimicking standard reductions at play in proofs of normalisation for Natural Deduction calculi \cite{prawitz1965}, compute detours connected to non-privileged arbitrary inferences, relative to premises obtained via privileged inferences. Roughly, an argument is valid when, by applying its reductions, it can be transformed (upon closure of its unbound assumptions) into a proof-structure where the privileged inferences play a prior role. Validity is referred to the very same atomic proof systems as in B-eS, and a formula being a consequence of a set of formulas (over a given proof system) is defined as existence of a valid (on the system) argument from the set of formulas as assumptions-set, to the formula as conclusion. This approach is normally called \emph{Proof-Theoretic Validity} (P-tV), and it can again be carried out in two fashions:

\begin{itemize}
    \item arguments may be valid in a \emph{monotonic} or \emph{non-monotonic} way, depending on whether validity on a given proof system is (resp. is not) preserved when extending the system, and
    \item arguments may be valid in an \emph{introduction-based} or \emph{elimination-based} way, depending on whether one takes the introduction (resp. elimination) rules as privileged.
\end{itemize}
Below, I shall focus only on (propositional) \emph{monotonic introduction-based} P-tV (miP-tV)---for approaches that privilege elimination rules (not necessarily in a monotonic context) see \cite{GheorghiuPym, hermogenespragmatist, sandqvist, takemura}.

Questions about completeness or incompleteness relative to miP-tV are normally dealt with by reading the miP-tV consequence relation in terms of the miB-eS one, and then transferring to miP-tV the results available for miB-eS---see e.g. \cite{piechaschroeder-heister2019}. It has been remarked, however, that this strategy might work to a limited extent only \cite{piccolominithesis, piccolominibook, piccolomininote, schroederheisterrolf}. For, due to certain structural differences between miB-eS and miP-tV, a reading of miP-tV consequence via miB-eS consequence may force a weakening of the kind of reductions at play in the miP-tV concept of argument.

Prawitz's original notion of reduction \cite{prawitz1973} seems to enjoy two main features: first, reductions must be \emph{constructive}, in some reasonable sense of constructivity for functions; second, reductions must be \emph{uniform}, in the rough sense that they must be described in such a way that their output values do not depend on atomic proof systems which their input values are valid on, relative to some assignment of given reduction-sets. If one gives these conditions up, then miB-eS can be proved to be equivalent to miP-tV in the sense that a formula is a consequence of a (finite) set of formulas in miB-eS (over an atomic proof system) if and only if the formula is a consequence of the set of formulas in miP-tV (over the same atomic proof system) \cite{piccolominiinversion, stafford1}.\footnote{The same is proved for the corresponding non-monotonic variants in \cite{piccolominisomeresults}, and for non-monotonic Prawitz's theory of grounds \cite{prawitz2015} in \cite{piccolominibook}, in both cases with a strategy which easily adapts to the monotonic framework.} Thus, the completeness or incompleteness results for miB-eS go through for miP-tV too, and in particular, $\texttt{IL}$ becomes incomplete over miP-tV \cite{piccolominiinversion}. But if the conditions at issue \emph{are not} given up, the equivalence between miB-eS and miP-tV consequence is harder to establish, and the (in)completeness of $\texttt{IL}$ over miP-tV remains unsettled.

In what follows, however, I single out a set of reductions such that (1) its elements are constructive and uniform (independently of how uniformity is more precisely defined), and (2) it validates the atomic Split rule, which is known to be underivable in $\texttt{IL}$. This is done by building upon a proof-strategy to be found in \cite{piechadecampossanz, piechadecampossanzschroederheisterpeirce, piechaschroeder-heisterdecampossanz} where, based on a notion of (intuitionistic) construction, a more Prawitz-friendly incompleteness result for $\texttt{IL}$ is established (namely, without translating directly from miB-eS). I shall also appeal to some insights into the notion of uniformity, in such a way as to cope with a proof by Pezlar of the constructive validity of a generalised Split rule (where the antecedent of the premise is a Harrop formula) \cite{pezlarselector}. Let me remark, in passing, that similar incompleteness results for non-monotonic P-tV are proved in \cite{piccolominiwem} and, in a more general fashion, in \cite{piccolominiprawitz}.

The paper is structured as follows. In Section 2, I define miB-eS, and sum up the above mentioned incompleteness theorem for $\texttt{IL}$ relative to it. In Section 3, I define miP-tV based on Prawitz's papers from the 70s. In Section 4, I sum up the main differences between miB-eS and miP-tV, with focus on the kind of problems one has when trying to read the notion of consequence of the latter via the notion of consequence at play in the former. In Sections 5 and 6, I discuss what completeness of $\texttt{IL}$ comes down to when formulated in strict miP-tV terms. Finally, in Section 7 I provide an miP-tV proof of the logical validity of the atomic Split rule, hence of incompleteness of $\texttt{IL}$ over miP-tV. In the concluding remarks I introduce some points for further discussion, especially concerning a notion of logical validity closed under substitutions of atoms with formulas.

\section{Base-extension Semantics and intuitionistic incompleteness}

Starting from miB-eS, I will limit myself, as said, to the propositional level. The language is hence defined thus.

\begin{definition}
    The \emph{propositional language} $\mathscr{L}$ is given by the following grammar:

    \begin{center}
        $A \coloneqq p_1, p_2, p_3, ... \ | \ \bot \ | \ A \wedge A \ | \ A \vee A \ | \ A \rightarrow A$
    \end{center}
    where $\bot$ is an atomic constant symbol for absurdity\footnote{For the purposes of this paper, it is immaterial whether $\bot$ is understood as an atom, as done here, or as a $0$-ary connective, as done e.g. in \cite{sandqvist}. The set of reductions I shall present in Section 7 which validate (logically) the atomic Split rule works in both cases. In general, however, the way how $\bot$ is dealt with, both linguistically and at a semantic level (see below), is crucial for completeness over PTS, especially miB-eS.}, and $\neg A \stackrel{\emph{def}}{=} A \rightarrow \bot$. The set $\emph{\texttt{FORM}}_\mathscr{L}$ of the \emph{formulas} of $\mathscr{L}$ contains the set $\emph{\texttt{ATOM}}_\mathscr{L} = \{p_i \ | \ i \in \mathbb{N}\} \cup \{\bot\}$ of the \emph{atoms} of $\mathscr{L}$.
\end{definition}
\noindent miB-eS (and miP-tV) require as said a preliminary notion of atomic proof system, which the notion of consequence (resp. of argumental validity) is relativised over. These proof systems are just sets of atomic rules of different levels of complexity.

\begin{definition}
    The notion of \emph{atomic rule of level} $n \geq 0$ is defined by induction as follows:
    \begin{itemize}
        \item with $A \in \emph{\texttt{ATOM}}_\mathscr{L}$,

        \begin{prooftree}
            \AxiomC{}
            \UnaryInfC{$A$}
        \end{prooftree}        
        is an atomic rule of level $0$
        \item if $A_1, ..., A_n, A \in \texttt{ATOM}_\mathscr{L}$, then

        \begin{prooftree}
            \AxiomC{$A_1$}
            \AxiomC{$\dots$}
            \AxiomC{$A_n$}
            \TrinaryInfC{$A$}
        \end{prooftree}
        is an atomic rule of level $1$
        \item if $A_1, ..., A_n, A \in \emph{\texttt{ATOM}}_\mathscr{L}$ and if $\Re_i$ are atomic rules whose maximal level is $k \ (i \leq n)$, then
        
        \begin{prooftree}
            \AxiomC{$[\Re_1]$}
            \noLine
            \UnaryInfC{$A_1$}
            \AxiomC{$\dots$}
            \AxiomC{$[\Re_n]$}
            \noLine
            \UnaryInfC{$A_n$}
            \TrinaryInfC{$A$}
        \end{prooftree}
        
        is an atomic rule of level $k + 2$.
    \end{itemize}
\end{definition}
\noindent As usual in a Natural Deduction setting, square brackets indicate binding of assumptions. However, what is discharged in the atomic rules above are not assumed formulas, but assumed atomic rules of lower level than the rule where the dischargement takes place. So, these are so-called \emph{higher-level rules} \cite{schroeder-heister1984}. Before going to the next definition, let me just introduce the set

\begin{prooftree}
    \AxiomC{$\bot$}
    \LeftLabel{$\texttt{AE} = \{$}
    \RightLabel{$ \ | \ \ A \in \texttt{ATOM}_\mathscr{L}\}$}
    \UnaryInfC{$A$}
\end{prooftree}

\begin{definition}
    An \emph{atomic base of level} $n \geq 0$ is a set of atomic rules $\mathfrak{B} = \{\mathfrak{r}_1, ..., \mathfrak{r}_n\} \cup \emph{\texttt{AE}}$ such that $n$ is maximal among the levels of $\mathfrak{r}_1, ..., \mathfrak{r}_n$.\footnote{Here, $n$ is taken to be a finite number, but it is immaterial for the purpose of this paper whether we take the atomic rules in a base to be finite or infinite in number---in other words, some atomic rules may also actually be rule \emph{schemes} (e.g., a rule of level 2 mimicking elimination of disjunction). Again, the incompleteness result of Section 7 abstracts from this issue.}
\end{definition}
\noindent I am assuming that atomic bases are, as one says, always of an intuitionistic kind \cite{piechaschroeder-heister2019}, i.e., that they always include an atomic explosion rule for each atom of $\mathscr{L}$.\footnote{In line with footnote 2 above, let me remark that the way one deals with $\bot$ over atomic bases is crucial for issues of completeness over PTS, especially miB-eS. In particular, rather than requiring atomic bases to be always of an intuitionistic kind, one could require that the validity of $\bot$ on the base is tantamount to the validity of every atom on the base. As was already the case with similar issues above, this is immaterial for the incompleteness result of Section 7. I shall not take into account here an approach recently discussed in \cite{barrosopereiranascimento}, where atomic bases (and their extensions) are always required to be consistent.} The atomic base $\mathfrak{B}^\emptyset = \texttt{AE}$ is also called the \emph{empty base}. Observe that the atomic explosion rules do not play any role when counting the complexity of the base. Hence, $\mathfrak{B}^\emptyset$ has complexity $0$.

\begin{definition}
    $\mathfrak{C}$ is an \emph{extension} of $\mathfrak{B}$ iff $\mathfrak{C} \supseteq \mathfrak{B}$.
\end{definition}

\begin{definition}
    The notion of \emph{atomic derivation} is defined by induction as follows:
    \begin{itemize}
        \item the single node derivation
        \begin{prooftree}
        \AxiomC{}
        \UnaryInfC{$A$}
    \end{prooftree}
is an atomic derivation
        \item if the following are atomic derivations
        \begin{prooftree}
            \AxiomC{$\mathfrak{G}_i, \Re_i$}
            \noLine
            \UnaryInfC{$\mathscr{D}_i$}
            \noLine
            \UnaryInfC{$A_i$}
        \end{prooftree}
        where $A_i$ is the premise of an atomic rule $\mathfrak{r}$ of the form
        \begin{prooftree}
            \AxiomC{$[\Re_1]$}
            \noLine
            \UnaryInfC{$A_1$}
            \AxiomC{$\dots$}
            \AxiomC{$[\Re_n]$}
            \noLine
            \UnaryInfC{$A_n$}
            \RightLabel{$\mathfrak{r}_2$}
            \TrinaryInfC{$A$}
        \end{prooftree}
        \emph{(}$i \leq n$\emph{)}, then
        \begin{prooftree}
            \AxiomC{$\mathfrak{G}_1, [\Re_1]$}
            \noLine
            \UnaryInfC{$\mathscr{D}_1$}
            \noLine
            \UnaryInfC{$A_1$}
            \AxiomC{$\dots$}
            \AxiomC{$\mathfrak{G}_n, [\Re_n]$}
            \noLine
            \UnaryInfC{$\mathscr{D}_n$}
            \noLine
            \UnaryInfC{$A_n$}
            \RightLabel{$\mathfrak{r}_2$}
            \TrinaryInfC{$A$}
    \end{prooftree}
        is an atomic derivation.
    \end{itemize}
\end{definition}

\begin{definition}
    $A$ is \emph{derivable} from $\mathfrak{G}$ in $\mathfrak{B}$, written $\mathfrak{G} \vdash_\mathfrak{B} A$, iff there is an atomic derivation $\mathscr{D}$ such that, for every $\mathfrak{r}$ applied and undischarged in $\mathscr{D}$, $\mathfrak{r} \notin \mathfrak{G} \Longrightarrow \mathfrak{r} \in \mathfrak{B}$. The set $\emph{\texttt{DER}}_\mathfrak{B}$ of the \emph{atomic derivations of} $\mathfrak{B}$ is the set of the atomic derivations $\mathscr{D}$ such that, for every $\mathfrak{r}$, if $\mathfrak{r}$ is applied and not discharged in $\mathscr{D}$, then $\mathfrak{r} \in \mathfrak{B}$.
\end{definition}
\noindent Observe that, if $A \in \mathfrak{B}$ we have both $A \vdash_\mathfrak{B} A$ and $\vdash_\mathfrak{B} A$. When $A \notin \mathfrak{B}$, we have instead the first result only. Similarly, if $\mathfrak{B} = \{p\}$, then $p, q, \mathfrak{r} \vdash_\mathfrak{B} r$ and $q, \mathfrak{r} \vdash_\mathfrak{B} r$, where $\mathfrak{r}$ is as in the derivation witnessing the result, i.e.,

\begin{prooftree}
    \AxiomC{}
    \UnaryInfC{$p$}
    \AxiomC{}
    \UnaryInfC{$q$}
    \RightLabel{$\mathfrak{r}$}
    \BinaryInfC{$r$}
\end{prooftree}
If $\mathfrak{B} = \{p, \mathfrak{r}\}$, besides the previous results, we have $q \vdash_\mathfrak{B} r$ too (by the same derivation as above). The above is an atomic derivation in $\mathfrak{B}$ when $\mathfrak{B} = \{p, q, \mathfrak{r}\}$. We can now define the variant of miB-eS that we shall focus on in this paper. Since I am restricting myself to the propositional level, I will allow myself to use quantifiers in the meta-language. In what follows, $\Gamma$ has to be understood as finite---this is to cope with some issues of compact monotonicity pointed out by \cite{stafford1} in connection with miP-tV.

\begin{definition}
    $A$ is a \emph{consequence} of $\Gamma \subset \emph{\texttt{FORM}}_\mathscr{L}$ \emph{over} $\mathfrak{B}$ \emph{in the miB-eS sense}, written $\Gamma \Vdash_\mathfrak{B} A$, iff
    \begin{itemize}
    \item $\Gamma = \emptyset \Longrightarrow$
    \begin{itemize}
        \item $A \in \emph{\texttt{ATOM}}_\mathscr{L} \Longrightarrow \ \vdash_\mathfrak{B} A$
        \item $A = B \wedge C \Longrightarrow \ \Vdash_\mathfrak{B} B$ and $\Vdash_\mathfrak{B} C$
        \item $A = B \vee C \Longrightarrow \ \Vdash_\mathfrak{B} B$ or $\Vdash_\mathfrak{B} C$
        \item $A = B \rightarrow C \Longrightarrow B \Vdash_\mathfrak{B} C$
    \end{itemize}
    \item $\Gamma \neq \emptyset \Longrightarrow \forall \mathfrak{C} \supseteq \mathfrak{B} \ (\Vdash_\mathfrak{C} \Gamma \Longrightarrow \ \Vdash_\mathfrak{B} A)$
    \end{itemize}
    where $\Vdash_\mathfrak{C} \Gamma$ means $\Vdash_\mathfrak{C} B$ for every $B \in \Gamma$.
\end{definition}

\begin{definition}
    $A$ is a \emph{logical consequence} of $\Gamma$ \emph{in the miB-eS sense}, written $\Gamma \Vdash A$, iff $\forall \mathfrak{B} \ (\Gamma \Vdash_\mathfrak{B} A)$.
\end{definition}
\noindent Observe that there is no special clause for $\bot$. This is essentially because we are assuming the atomic bases to be always intuitionistic. Remarkable facts about $\Vdash$ (also useful for what I shall say below) are the following---see \cite{piechadecampossanzschroederheisterpeirce, piechaschroeder-heisterdecampossanz, piechaschroeder-heister2019}.

\begin{proposition}
    $\Gamma \Vdash_\mathfrak{B} A$ iff $\forall \mathfrak{C} \supseteq \mathfrak{B} \ (\Gamma \Vdash_\mathfrak{C} A)$
\end{proposition}

\begin{proposition}
    $\Gamma \Vdash A$ iff $\Gamma \Vdash_{\mathfrak{B}^\emptyset} A$.
\end{proposition}

\begin{proposition}
    $\emph{\texttt{IL}}$ is sound over $\Vdash$, i.e., $\Gamma \vdash_\emph{\texttt{IL}} A \Longrightarrow \Gamma \Vdash A$.
\end{proposition}

\begin{definition}
    $\emph{\texttt{IL}}$ is \emph{complete} over $\Vdash$ iff $\Gamma \Vdash A \Longrightarrow \Gamma \vdash_{\emph{\texttt{IL}}} A$.
\end{definition}

\begin{theorem}
    $\emph{\texttt{IL}}$ is not complete over $\Vdash$, i.e., $\exists \ \Gamma, A \ (\Gamma \Vdash A$ and $\Gamma \not\vdash_{\emph{\texttt{IL}}} A)$.
\end{theorem}
\noindent The incompleteness results span over a number of papers, of which the most general one is \cite{piechaschroeder-heister2019}, but see also \cite{piechadecampossanz, piechadecampossanzschroederheisterpeirce, piechaschroeder-heisterdecampossanz}. It should be remarked, however, that de Campos Sanz, Piecha and Schroeder-Heister prove an actually stronger result, namely, that $\texttt{IL}$ is incomplete over $\Vdash^*$, where $\Vdash^*$ is $\Vdash$ closed under substitutions of atoms with formulas. Here, I stick to simple completeness, i.e., to Definition 9 and Theorem 1. I shall come back to the issue of closure under substitutions in the concluding remarks.

\section{Proof-Theoretic Validity}

Whereas in miB-eS one starts with the primitive consequence relation $\Vdash$, the consequence notion of miP-tV is a derived one, since it is defined in terms of existence of suitable valid arguments from given assumptions to a given conclusion. The miP-tV notion of valid argument is defined in this section for the propositional level again. And again, arguments will be either valid relative to an atomic base or, when validity obtains under universal quantification over atomic bases, they will be logically valid. Therefore, we keep Definitions 1, 2, 3, 4, 5 and 6 from Section 2. My presentation will be as close as possible to that of Prawitz's \cite{prawitz1973}, although with differences which I shall explain below.

An argument in Prawitz's original picture is a pair, where the first element, called argument structure, is a pair too. The left-hand element of this second pair is a derivational tree in the style of Gentzen's Natural Deduction, except that we put no restriction on the kind of inferences we are allowed to use (namely, inferences should not be conceived of as instantiating inference rules in a given Natural Deduction system). Thus, the nodes of the tree are labelled by formulas of $\mathscr{L}$, while the edges are to be understood as arbitrary inferential transitions. The formulas labelling the top-nodes are called the assumptions of the argument structure, while the root-node is called the conclusion of the argument structure. Some inferences may bind assumptions, which is why the argument structure pair involves a second element, that is, a function from top-nodes to lower nodes which indicates where in the argument structure the bindings of assumptions take place. Since Prawitz uses atomic bases of level $0$ or $1$ only, this is enough for defining argument structures in his original framework. However, we are now working with atomic bases of higher level, thus we need a bit more complex picture. The following definition does the job, while complying (hopefully) with Prawitz's original picture.

\begin{definition}
    An \emph{argument structure over $\mathscr{L}$} is a pair $\langle T, \langle f, h, g \rangle \rangle$ such that
    
    \begin{itemize} 
    \item $T$ is a finite rooted tree with order relation $\prec$, whose nodes are labelled by formulas of $\mathscr{L}$. The top-nodes of $T$ are partitioned into two groups
    \begin{itemize}
    \item assumption-labels, written $T^{\emph{\texttt{As}}}$
    \item axiom-labels, written $T^{\emph{\texttt{Ax}}}$
    \end{itemize}
    \item $f$ is a function defined on some $G \subseteq T^{\emph{\texttt{As}}}$ and is such that, $\forall \mu \in G$, it holds that $\mu \prec f(\mu)$
     \item $h$ is a function defined on some $H \subseteq T^{\emph{\texttt{ax}}}$ and is such that, $\forall \mu \in H$, $\mu$ is labelled by an atom, $\mu \prec h(\mu)$, $h(\mu)$ and all its children are labelled by atoms, and there is no $\nu \in G$ such that $f(\nu) = h(\mu)$;
    \item $g$ is a function defined on some $J \subseteq \mathcal{P}(\prec)$ such that, for every $K \in J$,
    \begin{itemize}
    \item $K$ contains all and only the edges that link a given node $\mu_K$ to all its children, and
    \item both $\mu_K$ and its children are labelled by atoms, and
    \item there is no $\nu \in G$ such that $f(\nu) = \mu_K$, and
    \item the function is such that $\forall K \in J$, it holds that $\mu_K \prec g(K)$, $g(K)$ and all its children are labelled by atoms, and there is no $\xi \in G$ such that $f(\xi) = g(K)$.
    \end{itemize}
    \end{itemize}
\end{definition}

\begin{definition}
    Given $\mathscr{D} = \langle T, \langle f, h, g \rangle \rangle$ with root-node labelled by $A$, the elements of $T^{\emph{\texttt{As}}}$ are the \emph{assumption-formulas} of $\mathscr{D}$, while $A$ is the \emph{conclusion} of $\mathscr{D}$.
\end{definition}

\noindent The meaning of $T^{\texttt{Ax}}$ in Definition 10 is that the formulas of the given top-nodes in $T$ are not assumptions, but axioms. Graphically, they are represented thus

\begin{prooftree}
    \AxiomC{}
    \UnaryInfC{$A$}
\end{prooftree}
where $A$ is meant to depend on no assumptions. It is worth remarking that, in virtue of Definition 5, top-nodes in atomic derivations are \emph{always} axioms (i.e., atomic rules of level $0$). The meaning of $f, h, g$ in the same definition is, in Natural Deduction terminology, that they are discharge functions---see \cite{prawitz1965, schroeder-heister1984, schroeder-heister2006}. Top-nodes $\mu \in T^{\texttt{as}}$ in the domain of $f$ are assumptions discharged in $T$, while top-nodes $\mu \in T^{\texttt{ax}}$ and sets of edges $K \in J$ in the domain of $h, g$ respectively are assumed atomic rules (of level $0$ with $h$) discharged by an atomic rule. In all cases, the dischargement takes place at the nodes $f(\mu), h(\mu), g
(K)$.

\begin{definition}
    $\mathscr{D}$ is \emph{closed} iff all its assumption-formulas are discharged, and it is \emph{open} otherwise.
\end{definition}

\begin{definition}
 Where $\Gamma$ is the set of the undischarged assumption-formulas of $\mathscr{D}$ and $A$ is the conclusion of $\mathscr{D}$, $\mathscr{D}$ is an argument structure \emph{from $\Gamma$ to} \emph{(}or \emph{for)} $A$.
\end{definition}

\noindent  I shall assume as defined (or definable) in a suitably precise way the notions of \emph{(immediate) sub-structure of $\mathscr{D}$}, and of \emph{substitution of the sub-structure $\mathscr{D}^*$ with the structure $\mathscr{D}^{**}$ in $\mathscr{D}$}---written $\mathscr{D}[\mathscr{D}^{**}/\mathscr{D}^*]$. On the one hand, the meaning of these notions is obvious. On the other hand, one must in both cases take care about re-definition of discharge functions. E.g., concerning the second notion, since $\mathscr{D}[\mathscr{D}^{**}/\mathscr{D}^*]$ may not be a sub-structure of $\mathscr{D}$, the discharge functions from $\mathscr{D}$ may not be ‘‘transferable" to $\mathscr{D}[\mathscr{D}^{**}/\mathscr{D}^*]$. The latter must nonetheless be such that formulas or atomic rules discharged at a node $\mu$ in $\mathscr{D}$, are discharged in $\mathscr{D}[\mathscr{D}^{**}/\mathscr{D}^*]$---if any---at a node $\mu^*$ that $\mu$ is ‘‘mapped onto" in $\mathscr{D}[\mathscr{D}^{**}/\mathscr{D}^*]$, whereas some dischargements in $\mathscr{D}$ may ‘‘disappear" in $\mathscr{D}[\mathscr{D}^{**}/\mathscr{D}^*]$. I will allow myself to abstract from these details in what follows, but see \cite{prawitz1965, schroeder-heister1984, schroeder-heister2006} for more.

\begin{definition}
    Let $\mathscr{D}$ from $\Gamma = \{A_1, ..., A_n\}$ to $A$, and $\sigma$ a function from and to argument structures such that $\sigma(A_i)$ is a \emph{(}closed\emph{)} argument structure for $A_i$ \emph{(}$i \leq n$\emph{)}. Then, $\mathscr{D}^\sigma = \mathscr{D}[\sigma(A_1), ..., \sigma(A_n)/A_1, ..., A_n]$ is called the \emph{(closed) $\sigma$-instance of $\mathscr{D}$}.
\end{definition}

\begin{definition}
    An \emph{inference} is a triple $\langle \langle \mathscr{D}_1, ..., \mathscr{D}_n \rangle, A, \delta \rangle$, where $\delta$ is an extension of the discharge functions associated to the $\mathscr{D}_i$-s ($i \leq n$). The \emph{argument structure associated to the inference}, indicated by the figure
    \begin{prooftree}
        \AxiomC{$\mathscr{D}_1$}
        \AxiomC{$\dots$}
        \AxiomC{$\mathscr{D}_n$}
        \RightLabel{$\delta$}
        \TrinaryInfC{$A$}
    \end{prooftree}
    is obtained by conjoining the trees of $\mathscr{D}_i$-s ($i \leq n$) through a root-node $A$, and by extending via $\delta$ the discharge functions of the $\mathscr{D}_i$-s ($i \leq n$).
\end{definition}
\noindent The extension of the discharge functions is to be understood as underlying the same technical precautions as mentioned above, in connection with the notions of (immediate) sub-structure and of replacement of sub-structures with structures. Again, I shall abstract from these details, but see \cite{prawitz1965, schroeder-heister1984, schroeder-heister2006} for more.

\begin{definition}
    An \emph{inference rule} is a set of inferences, whose elements are called \emph{instances} of the rule.
\end{definition}

\noindent E.g., introduction (abstracting from the indexing of occurrences of assumptions) and elimination of implication are

\begin{center}
$\rightarrow_I \ = \ \{\langle \mathscr{D}, A \rightarrow B, f(A) = A \rightarrow B \rangle \ | \ \mathscr{D} \ \text{for} \ B\}$
\end{center}

\begin{center} 
$\rightarrow_E \ = \ \{\langle \mathscr{D}_1, \mathscr{D}_2, B, \emptyset \rangle \ | \ \mathscr{D}_1 \ \text{for} \ A \rightarrow B \ \text{and} \ \mathscr{D}_2 \ \text{for} \ A\}$
\end{center}
Each instance of these has associated argument structures
\begin{prooftree}
    \AxiomC{$[A]_1$}
    \noLine
    \UnaryInfC{$\mathscr{D}$}
    \noLine
    \UnaryInfC{$B$}
    \RightLabel{$1$}
    \UnaryInfC{$A \rightarrow B$}
    \AxiomC{$\mathscr{D}_1$}
    \noLine
    \UnaryInfC{$A \rightarrow B$}
    \AxiomC{$\mathscr{D}_2$}
    \noLine
    \UnaryInfC{$A$}
    \BinaryInfC{$B$}
    \noLine
    \BinaryInfC{}
\end{prooftree}
respectively. Thus, the rules themselves can be given the respective, usual meta-linguistic descriptions\footnote{A similar notation is used below also for referring to schematic argument structures (i.e., argument structures whose formulas share a common logical form), and in some cases even for non-schematic argument structures (i.e., argument structures involving specific formulas). The context will make clear whether I am speaking of rules or of argument structures---where it is to be observed that, of course, a one-step instance of a rule which does not discharge assumptions \emph{is} an argument structure.}

\begin{prooftree}
    \AxiomC{$[A]$}
    \noLine
    \UnaryInfC{$B$}
    \UnaryInfC{$A \rightarrow B$}
    \AxiomC{$A \rightarrow B$}
    \AxiomC{$A$}
    \BinaryInfC{$B$}
    \noLine
    \BinaryInfC{}
\end{prooftree}
\noindent In what follows, I shall harmlessly identify a rule with the set of the argument structures associated to its instances.

\begin{definition}
    $\mathscr{D}$ is \emph{canonical} iff it is associated to an instance of an introduction rule, and it is \emph{non-canonical} otherwise.
\end{definition}

\begin{definition}
    Given an inference rule $\mathscr{R}$, a \emph{reduction} for $\mathscr{R}$ is a constructive function $\phi$ such that:
    \begin{itemize}
        \item[1.] $\phi$ goes from argument structures to argument structures
        \item[2.] $\phi$ is defined on some $\mathbb{D} \subseteq \mathscr{R}$ such that $\forall \sigma$, $\mathscr{D} \in \mathbb{D} \Longrightarrow \mathscr{D}^\sigma \in \mathbb{D}$
        \item[3.] $\forall \mathscr{D} \in \mathbb{D}$, $\mathscr{D}$ is from $\Gamma$ to $A \Longrightarrow \phi(\mathscr{D})$ is from $\Gamma^* \subseteq \Gamma$ to $A$
        \item[4.] $\forall \sigma$, $\phi(\mathscr{D}^\sigma) = \phi(\mathscr{D})^\sigma$
    \end{itemize}
\end{definition}
\noindent So for example, the standard reduction $\phi_\rightarrow$ for elimination of implication is defined on

\begin{center}
    $\{\langle \mathscr{D}_1, \mathscr{D}_2, B, \emptyset \rangle \ | \ \mathscr{D}_1 \in \ \rightarrow_I\} \subset \ \rightarrow_E$
\end{center}
such that, as usual,

\begin{prooftree}
    \AxiomC{$[A]_1$}
    \noLine
    \UnaryInfC{$\mathscr{D}^*_1$}
    \noLine
    \UnaryInfC{$B$}
    \RightLabel{$1$}
    \UnaryInfC{$A \rightarrow B$}
    \AxiomC{$\mathscr{D}_2$}
    \noLine
    \UnaryInfC{$A$}
    \BinaryInfC{$B$}
    \AxiomC{$\stackrel{\phi_\rightarrow}{\Longrightarrow}$}
    \noLine
    \UnaryInfC{}
    \AxiomC{$\mathscr{D}_2$}
    \noLine
    \UnaryInfC{$[A]$}
    \noLine
    \UnaryInfC{$\mathscr{D}^*_1$}
    \noLine
    \UnaryInfC{$B$}
    \noLine
    \TrinaryInfC{}
\end{prooftree}

\begin{definition}
    With $\mathfrak{J}$ set of reductions, an \emph{extension} of $\mathfrak{J}$ is a set of reductions $\mathfrak{H}$ such that $\mathfrak{H} \supseteq \mathfrak{J}$.
\end{definition}

\begin{definition}
    $\mathscr{D}$ \emph{immediately reduces} to $\mathscr{D}^*$ \emph{relative to} $\mathfrak{J}$ iff, for some sub-structure $\mathscr{D}^{**}$ of $\mathscr{D}$, there is $\phi \in \mathfrak{J}$ such that $\mathscr{D}^* = \mathscr{D}[\phi(\mathscr{D}^{**})/\mathscr{D}^{**}]$. The relation of $\mathscr{D}$ \emph{reducing} to $\mathscr{D}^*$ \emph{relative to} $\mathfrak{J}$, written $\mathscr{D} \leq_\mathfrak{J} \mathscr{D}^*$, is the reflexive-transitive closure of the relation of immediate reducibility of $\mathscr{D}$ to $\mathscr{D}^*$.
\end{definition}

\begin{definition}
    $\langle \mathscr{D}, \mathfrak{J} \rangle$ is \emph{valid} on $\mathfrak{B}$ iff
    \begin{itemize}
        \item $\mathscr{D}$ is closed $\Longrightarrow$
        \begin{itemize}
            \item the conclusion of $\mathscr{D}$ is an atom $A$ $\Longrightarrow \mathscr{D} \leq_\mathfrak{J} \mathscr{D}^*$ with $\mathscr{D}^*$ witnessing $\vdash_\mathfrak{B} A$;
            \item the conclusion of $\mathscr{D}$ is not atomic $\Longrightarrow \mathscr{D} \leq_\mathfrak{J} \mathscr{D}^*$ for $\mathscr{D}^*$ canonical with immediate sub-structures valid on $\mathfrak{B}$ when paired with $\mathfrak{J}$;
        \end{itemize}
        \item $\mathscr{D}$ is open $\Longrightarrow \forall \sigma, \forall A_1, ..., A_n \in \Gamma, \forall \mathfrak{H} \supseteq \mathfrak{J}$ and $\forall \mathfrak{C} \supseteq \mathfrak{B}$, if $\langle \sigma(A_i), \mathfrak{H} \rangle$ is valid on $\mathfrak{C} \ \emph{(}i \leq n\emph{)}$, then $\langle \mathscr{D}^\sigma, \mathfrak{H} \rangle$ is valid on $\mathfrak{C}$.
    \end{itemize}
\end{definition}

\begin{definition}
    $A$ is a \emph{consequence} of $\Gamma$ \emph{over} $\mathfrak{B}$ \emph{in the miP-tV sense}, written $\Gamma \models_\mathfrak{B} A$, iff $\exists \ \langle \mathscr{D}, \mathfrak{J} \rangle$ valid on $\mathfrak{B}$ with $\mathscr{D}$ from $\Gamma$ to $A$.
\end{definition}

\begin{definition}
    $\langle \mathscr{D}, \mathfrak{J} \rangle$ is \emph{logically valid} iff it is valid on every $\mathfrak{B}$.
\end{definition}

\begin{definition}
     $A$ is a \emph{logical consequence} of $\Gamma$ \emph{in the miP-tV sense}, written $\Gamma \models A$, iff $\exists \ \langle \mathscr{D}, \mathfrak{J} \rangle$ logically valid with $\mathscr{D}$ from $\Gamma$ to $A$.
\end{definition}

Results similar to those for $\Vdash$ in Propositions 1, 2 and 3, can be proved for the notion of (logically) valid argument and for $\models_{(\mathfrak{B})}$ too.\footnote{I use the notation $\models_{(\mathfrak{B})}$ [resp. $\Vdash_{(\mathfrak{B})}$] when I mean to refer what I am saying indifferently to either $\models_\mathfrak{B}$ [resp. $\Vdash_\mathfrak{B}$], for some $\mathfrak{B}$, or to $\models$ [resp. $\Vdash$].} I.e.: the approach is monotonic; logical validity and logical consequence are, respectively, tantamount to validity and consequence over the empty base; $\texttt{IL}$ is sound. Proofs are omitted but see \cite{piccolominiinversion} and \cite{schroeder-heister2006}---although the latter uses a different notion of reduction.

\begin{proposition}
    $\langle \mathscr{D}, \mathfrak{J} \rangle$ valid on $\mathfrak{B} \Longleftrightarrow \forall \mathfrak{C} \supseteq \mathfrak{B}, \langle \mathscr{D}, \mathfrak{J} \rangle$ valid on $\mathfrak{C}$.
\end{proposition}

\begin{corollary}
    $\Gamma \models_\mathfrak{B} A \Longleftrightarrow \forall \mathfrak{C} \supseteq \mathfrak{B}, \Gamma \models_\mathfrak{C} A$.
\end{corollary}

\begin{proposition}
    $\langle \mathscr{D}, \mathfrak{J} \rangle$ logically valid $\Longleftrightarrow \langle \mathscr{D}, \mathfrak{J} \rangle$ valid on $\mathfrak{B}^\emptyset$.
\end{proposition}

\begin{corollary}
    $\Gamma \models A \Longleftrightarrow \Gamma \models_{\mathfrak{B}^\emptyset} A$.
\end{corollary}

\begin{proposition}
    $\emph{\texttt{IL}}$ is sound over $\models$, i.e., $\Gamma \vdash_{\emph{\texttt{IL}}} A \Longrightarrow \Gamma \models A$.
\end{proposition}

\noindent Another trivial result which I shall use below is monotonicity over sets of reductions---though there is no counter-part of it at the level of consequence, due to the absence of references to sets of reductions in this notion.

\begin{proposition}
    $\langle \mathscr{D}, \mathfrak{J} \rangle$ valid on $\mathfrak{B}$ \emph{[}resp. logically valid\emph{]} $\Longleftrightarrow \forall \mathfrak{H} \supseteq \mathfrak{J}, \ \langle \mathscr{D}, \mathfrak{H} \rangle$ valid on $\mathfrak{B}$ \emph{[}resp. logically valid\emph{]}.
\end{proposition}

The presentation of miP-tV I gave here differs in essentially three respects from the original framework of Prawitz \cite{prawitz1973}. First, as said, Prawitz works with atomic bases of level at most $1$, while I assume that they contain atomic rules of higher-level. Second, Prawitz assumes atomic bases to be always consistent, whereas I am taking them to be always of an intuitionistic kind---i.e., to contain all atomic explosion rules. Third, Prawitz's version is \emph{non-monotonic}, in that the clause for the case when $\mathscr{D}$ is open in Definition 21 is given by him without extensions of $\mathfrak{B}$, so that, e.g., the one-step argument-structure from $p$ to $q$ is valid on $\mathfrak{B}^\emptyset$ with respect to the empty set of reductions, but fails to be so over $\mathfrak{B}^\emptyset \cup \{p\}$. This means that, in Prawitz's original framework, Proposition 4, Corollary 1, Proposition 5, and Corollary 2 fail. As I shall remark again in due course, the results from Section 7 hold also when one restricts oneself to bases of level $\leq 1$ of the kind described above.\footnote{And independently of how one deals with $\bot$ in the atomic bases, see footnotes 2 and 4 above.} The third difference is instead much more substantial. But I shall leave the non-monotonic approach aside here---however, see \cite{piccolomininote, piccolominisomeresults, piccolominiwem, piccolominiprawitz, piechaschroeder-heister2019, schroederheisterrolf}.

\section{Constructivity, uniformity, and admissibility-clause}

As said, $\texttt{IL}$ is incomplete over $\Vdash$ (Theorem 1). The question may now arise whether the incompleteness results provable relative to $\Vdash$ can be somehow ‘‘transferred" to $\models$, which would obtain if we could prove some result to the effect that $\Gamma \Vdash A \Longleftrightarrow \Gamma \models A$. However, this might not come without substantial limitations. First, we must limit ourselves to $\Gamma$ finite since, as pointed out by \cite{stafford1}, $\models_{(\mathfrak{B})}$ is only compactly monotone. Another critical point has instead to do with the sets of reductions that we use for justifying structures in valid arguments which witness $\models_{(\mathfrak{B})}$.

As stated in \cite{piccolominiinversion, stafford1}, $\Gamma$ always finite plus ‘‘liberal" sets of reductions yields $\Gamma \Vdash_\mathfrak{B} A \Longleftrightarrow \Gamma \models_\mathfrak{B} A$. For, under these conditions, one can show that $\models_\mathfrak{B}$ satisfies the same clauses as those which define $\Vdash_\mathfrak{B}$ in Definition 7. Via Proposition 2 and Corollary 2, we can then infer $\Gamma \Vdash A \Longleftrightarrow \Gamma \models A$, which is enough for an miP-tV variant of Theorem 1, i.e., for incompleteness of $\texttt{IL}$ over $\models$ \cite{piccolominiinversion}.

When referred to miP-tV, the constraint on $\Gamma$ finite is seemingly not worrying, since $\models_{(\mathfrak{B})}$ is defined as existence of a valid (on $\mathfrak{B}$) $\langle \mathscr{D}, \mathfrak{J} \rangle$, where $\mathscr{D}$ contains a finite number of assumptions only. To see whether ‘‘liberal" sets of reductions are equally acceptable or not, instead, it is time I explain what I mean more precisely by ‘‘liberal".

In the proof of $\Gamma \Vdash_\mathfrak{B} A \Longleftrightarrow \Gamma \models_\mathfrak{B} A$ referred to above, the ‘‘liberality" required for sets of reductions is most visible when showing that $\models_\mathfrak{B}$ satisfies the same condition as $\Vdash_\mathfrak{B}$ when $\Gamma \neq \emptyset$, i.e.,

\begin{itemize}
    \item[(*)] $\Gamma \models_\mathfrak{B} A \Longleftrightarrow \forall \mathfrak{C} \supseteq \mathfrak{B} \ (\models_\mathfrak{C} \Gamma \Longrightarrow \ \models_\mathfrak{C} A)$, with $\Gamma \neq \emptyset$.
\end{itemize}
While the left-to-right direction of (*) is quite trivial, the right-to-left direction requires associating to

\begin{prooftree}
    \AxiomC{$A_1$}
    \AxiomC{$\dots$}
    \AxiomC{$A_n$}
    \TrinaryInfC{$A$}
\end{prooftree}
---with $\Gamma = \{A_1, ..., A_n\}$---the set $\Pi$ of ‘‘pointers" $\phi$ mapping closed $\langle \mathscr{D}_i, \mathfrak{J}_i \rangle$ for $A_i$ valid on some $\mathfrak{C} \supseteq \mathfrak{B}$ onto a closed $\langle \mathscr{D}, \mathfrak{J} \rangle$ for $A$ valid on $\mathfrak{C}$---whose existence is guaranteed by the assumption $\forall \mathfrak{C} \supseteq \mathfrak{B} \ (\models_\mathfrak{C} \Gamma \Longrightarrow \ \models_\mathfrak{C} A)$, and whose form eventually depends on each of the $\langle \mathscr{D}_i, \mathfrak{J}_i \rangle$-s $(i \leq n)$.

More precisely, let $\mathbb{B}$ be the set of the extensions of an atomic base $\mathfrak{B}$. With $\mathfrak{C} \in \mathbb{B}$, let $\mathbb{J}_\mathfrak{C}$ be the set of the sequences $r = \langle \mathfrak{J}_1, ..., \mathfrak{J}_n \rangle$ of the sets of reductions such that, for some sequence of argument structures $d = \langle \mathscr{D}_1, ..., \mathscr{D}_n \rangle$, the $i$-th element of $d$ is for $A_i$, and is valid on $\mathfrak{C}$ when paired with the $i$-th element of $r$. We also say that $d$ \emph{matches} $r$ \emph{on} $\mathfrak{C}$ and, for any $r \in \mathbb{J}_\mathfrak{C}$, we consider the set $\mathbb{D}^r_\mathfrak{C}$ of the $d$-s that match $r$ on $\mathfrak{C}$. Because of the assumption $\forall \mathfrak{C} \supseteq \mathfrak{B} \ (\models_\mathfrak{C} \Gamma \Longrightarrow \ \models_\mathfrak{C} A)$, for every $\mathfrak{C} \in \mathbb{B}$, every $r \in \mathbb{J}_\mathfrak{C}$ and every $d \in \mathbb{D}^r_\mathfrak{C}$, there is $\langle \mathscr{D}^{r, d}, \mathfrak{J}^{r, d} \rangle$ for $A$ valid on $\mathfrak{C}$. So, for every $\mathfrak{C} \in \mathbb{B}$, every $r \in \mathbb{J}$ and every $d \in \mathbb{D}^r_\mathfrak{C}$, we define a function $\phi^r_\mathfrak{C}$ from $\mathbb{D}^r_\mathfrak{C}$ to the set of the argument structures by putting $\phi^r_\mathfrak{C}(d) = \mathscr{D}^{r, d}$, and a function $\chi^r_\mathfrak{C}$ from $\mathbb{D}^r_\mathfrak{C}$ to the set of the sets of reductions by putting $\chi^r_\mathfrak{C}(d) = \mathfrak{J}^{r, d}$. We now put

\begin{center}
    $\Pi = \bigcup_{\mathfrak{C} \in \mathbb{B}, r \in \mathbb{J}_\mathfrak{C}} (\phi^r_\mathfrak{C} \cup  \bigcup (\chi^r_\mathfrak{C}(\mathbb{D}^r_\mathfrak{C}))$.
\end{center}
---namely, the union set of all the $\phi^r_\mathfrak{C}$-s and of all the union sets of the images of all the $\mathbb{D}^r_\mathfrak{C}$-s via all the $\chi^r_\mathfrak{C}$-s, with $\mathfrak{C} \in \mathbb{B}$ and $r \in \mathbb{J}_\mathfrak{C}$. All the $\phi^r_\mathfrak{C}$-s can be easily shown to respect the conditions of definition 21, so they are reductions; likewise, the elements of the union set of the images of all the $\mathbb{D}^r_\mathfrak{C}$-s via all the $\chi^r_\mathfrak{C}$-s are reductions by assumption. Next, suppose we are given any closed $\langle \mathscr{D}_i, \mathfrak{J}_i \rangle$ for $A_i$ valid on any $\mathfrak{C} \supseteq \mathfrak{B}$, so $r = \langle \mathfrak{J}_1, ..., \mathfrak{J}_n \rangle \in \mathbb{J}_\mathfrak{C}$ and $d = \langle \mathscr{D}_1, ..., \mathscr{D}_n \rangle \in \mathbb{D}^r_\mathfrak{C}$. Since we have associated $\Pi$ to the one-step rule from $A_1, ..., A_n$ to $A$, the argument structure 

\begin{prooftree}
    \AxiomC{$\mathscr{D}_1$}
    \noLine
    \UnaryInfC{$A_1$}
    \AxiomC{$\dots$}
    \AxiomC{$\mathscr{D}_n$}
    \noLine
    \UnaryInfC{$A_n$}
    \TrinaryInfC{$A$}
\end{prooftree}
will reduce via $\Pi$ to the closed argument structure $\phi^r_\mathfrak{C}(\langle \mathscr{D}_1, ..., \mathscr{D}_n \rangle)$ for $A$, valid on $\mathfrak{C}$ when paired with $\chi^r_\mathfrak{C}(\langle \mathscr{D}_1, ..., \mathscr{D}_n \rangle)$, where $\phi^r_\mathfrak{C} \in \Pi$ and $\phi \in \Pi$, for every $\phi \in \chi^r_\mathfrak{C}(\langle \mathscr{D}_1, ..., \mathscr{D}_n \rangle)$.\footnote{Clearly, an easier solution would be to let $\Pi$ be simply the set of all the mappings from and to argument structures which respect the conditions of Definition 21. This surely contains the mappings $\phi^r_\mathfrak{C}$-s mentioned above, and all the elements of the sets of reductions with respect to which closed argument structures for $A$ are valid on $\mathfrak{C}$, for any $\mathfrak{C} \supseteq \mathfrak{B}$. This is the strategy adopted in \cite{piccolominiinversion}. The fact that all the $\phi^r_\mathfrak{C}$-s belong to the set of all the reductions stems from an (implicit) constructive reading of the meta-logical constants in $\forall \mathfrak{C} \supseteq \mathfrak{B} \ (\models_\mathfrak{C} \Gamma \Longrightarrow \ \models_\mathfrak{C} A)$, which is also used above. The idea is that, for any $\mathfrak{C} \supseteq \mathfrak{B}$, the reading yields an effective function from (sequences of) pairs whose first element is an argument structure, and whose second element is a proof of the fact that that argument structure is valid on $\mathfrak{C}$, to similar pairs. If one assumes that meta-proofs are unique, one can extract from this a function from and to the first elements of these pairs, with the required properties. With classical logic in the meta-language, the proof is easier: for all $\mathfrak{C} \supseteq \mathfrak{B}$, if there $i \leq n$ such that $\not\Vdash_\mathfrak{C} A_i$, then choose for $\mathfrak{C}$ the empty function; otherwise, pick a closed $\langle \mathscr{D}, \mathfrak{J} \rangle$ for $A$ valid on $\mathfrak{C}$, and define a constant function $\kappa_\mathfrak{C}$ onto $\mathscr{D}$, which maps onto $\mathscr{D}$ the first elements of all the sequences of closed $\langle \mathscr{D}_i, \mathfrak{J}_i \rangle$ for $A_i$ valid on $\mathfrak{C}$.} If one now compares the given $\Pi$ with, say, $\{\phi_\rightarrow, \iota\}$, where $\iota$ is

\begin{prooftree}
    \AxiomC{$\mathscr{D}$}
    \noLine
    \UnaryInfC{$A \rightarrow (B \rightarrow C)$}
    \UnaryInfC{$B \rightarrow (A \rightarrow B)$}
    \AxiomC{$\stackrel{\iota}{\Longrightarrow}$}
    \noLine
    \UnaryInfC{}
    \AxiomC{$\mathscr{D}$}
    \noLine
    \UnaryInfC{$A \rightarrow (B \rightarrow C)$}
    \AxiomC{$[A]_1$}
    \BinaryInfC{$B \rightarrow C$}
    \AxiomC{$[B]_2$}
    \BinaryInfC{$C$}
    \RightLabel{$1$}
    \UnaryInfC{$A \rightarrow C$}
    \RightLabel{$2$}
    \UnaryInfC{$B \rightarrow (A \rightarrow C)$}
    \noLine
    \TrinaryInfC{}
\end{prooftree}
(and where $\phi_\rightarrow$ is the standard reduction for elimination of implication defined above) one realises that $\Pi$ is in a sense much less schematic than $\{\phi_\rightarrow, \iota\}$. First, $\Pi$ might well be infinite (there might be as many ‘‘pointers" as are the closed valid arguments for the elements of $\Gamma$), while $\{\phi_\rightarrow, \iota\}$ is finite. Now, the fact that $\Pi$ may be infinite may not be a problem in itself---in fact, in \cite{prawitz1973} Prawitz does not officially require sets of reductions to be finite. The point is, though, that the kind of infinity we are dealing with here seems also to bring some ingredient of non-\emph{constructivity} in. 

Constructivity is something that in \cite{prawitz1973} Prawitz  requires for \emph{reductions} (as opposed to sets of reductions), the idea being that reductions have to be effectively computable functions in the same sense as $\phi_\rightarrow$ and $\iota$ are. The reason behind this, anyway, seems to be that it must always be possible to effectively carry out the reduction sequences which a given $\mathfrak{J}$ induces on a given (closed instance of) $\mathscr{D}$, whenever the functions in $\mathfrak{J}$ work as reductions for the rules instantiated in (the closed instance of) $\mathscr{D}$. So, it seems reasonable to require that the sets of reductions which $\models_{(\mathfrak{B})}$ is defined over, be constructive in the sense of allowing each time for an effective computation of the reduction sequences they generate.

This is easily seen to happen with $\{\phi_\rightarrow, \iota\}$, but what about $\Pi$? Is it true that, given any closed instance of the one-step open argument structure from $A_1, ..., A_n$ to $A$ above, we are effectively able to compute its value via $\Pi$? This seems to be hardly the case. Of course, $\Pi$ as a whole is not defined on \emph{every} closed instance of the one-step structure, since each ‘‘pointer" is defined only on closed $\langle \mathscr{D}_i, \mathfrak{J}_i \rangle$ for $A_i \ (i \leq n)$ valid on some $\mathfrak{C} \supseteq \mathfrak{B}$. But this is not worrying since, as per Definition 18, we are not constrained to require the domain of a reduction for a set of argument structures to coincide with such set---as is already the case with, e.g., $\phi_\rightarrow$.\footnote{In other words, since a rule is a set of inferences, and since an inference is just a triple whose first element is a sequence of argument structures, whose second element is a formula, and whose third element is a discharge function, by Definition 18 we may well take each ‘‘pointer" to be defined on the sub-set of the one-step rule from $A_1, ..., A_n$ to $A$ which contains just the inferences $\langle \langle \mathscr{D}_1, ..., \mathscr{D}_n \rangle, A, \emptyset \rangle$, with $\mathscr{D}_i$ closed for $A_i$, and valid relative to some $\mathfrak{J}_i$ on some $\mathfrak{C} \supseteq \mathfrak{B} \ (i \leq n)$. If we stick to Prawitz's definition, this is perfectly legitimate, and it implies that the elements of $\Pi$ do not operate on the whole set of the argument structures which end by applying the inference from $A_1, ..., A_n$ to $A$.} Nor can we say that individual ‘‘pointers" are not constructive, for a ‘‘pointer" is a constant function, and constant functions count as constructive. The trouble is that we are not provided with enough information for choosing the right ‘‘pointer", besides being granted that one such ‘‘pointer" surely belongs to our infinite set of reductions. So, while each element of $\Pi$ is constructive, as required in \cite{prawitz1973}, $\Pi$ itself seems not to be constructive in the above sense, i.e., it does not allow for an effective computation of closed instances of the argument structure it is associated to.

Non-constructivity is not the only problem with $\Pi$. Another issue is that, in the terminology of \cite{schroederheisterrolf}, $\Pi$ may be also said to be \emph{non-uniform}.\footnote{In \cite{piccolomininote}, I raised this issue relative to a Prawitzian interpretation of non-monotonic introduction-based B-eS (niB-es)---‘‘Prawitzian" meaning that, as I do below for $\Vdash_{(\mathfrak{B})}$, I interpreted the primitive notion of consequence at play in niB-eS as existence of suitable argument structures and reductions. In the terminology I employ there, the sets of reductions were said to be \emph{non-schematic}, but I think that Schroeder-Heister's terminological choice, in terms of non-uniformity, is better. The point I am discussing here for $\Vdash_{(\mathfrak{B})}$ is in any case the same as in \cite{piccolomininote}, the only difference being that, due to non-monotonicity, in niB-eS the problem lurks out when going from validity over an atomic base to logical validity---via the quantifiers-inversion issue I pinpoint below.} A set of reductions $\mathfrak{J}$ counts as uniform if, morally, it ‘‘works in the same way" over all $\mathfrak{B}$-s. More precisely, each $\phi \in \mathfrak{J}$ must be defined in such a way that the value of $\phi$ for any $\mathscr{D}$ in its domain, can be effectively indicated without references to structures which $\mathscr{D}$ reduces to on some specific $\mathfrak{B}$ relative to some other specific set of reductions, or to the validity of $\mathscr{D}$ or $\phi(\mathscr{D})$ under the same conditions.\footnote{Uniformity might be more strongly understood in the sense that the value of a reduction must not depend on \emph{further sets of reductions} relative to which input values compute or are valid, over given atomic bases. While it is fairly clear that, in \cite{prawitz1973}, Prawitz has some uniformity constraint as \emph{base-independence} in the background, it is not clear whether this extends to \emph{independence from sets of reductions}, thus I will leave this issue open in this paper. Be that as it may, the reductions I define in Section 7, for validating the atomic Split rule, are uniform in both senses. See also what Schroeder-Heister calls the requirement of \emph{full specification} in \cite{schroeder-heisterbasicideas}, a point which I shall also discuss below.} To put it with Schroeder-Heister,

\begin{quote}
    proof semantics can be formulated in ways that cannot be captured by sentence semantics [...] in the [...] proof semantics in Prawitz's \cite{prawitz1973} and \cite{prawitz1974}, the transformation of validity is something that is required to happen in a \emph{uniform} way. It is demanded that there is a single reduction procedure which does not depend on specific atomic systems but works for all [$\mathfrak{B}$-s] in the same way, i.e., schematically. \cite[503]{schroederheisterrolf}
\end{quote}
Again, it is easy to see that $\{\phi_\rightarrow, \iota\}$ is uniform in the sense required, but what about $\Pi$? Does it operate uniformly over all $\mathfrak{B}$-s? Also in this case, the answer seems to be negative. The way in which the ‘‘pointers" generate their values out of the closed $\langle \mathscr{D}_i, \mathfrak{J}_i \rangle$ for $A_i$ valid on some $\mathfrak{C} \supseteq \mathfrak{B}$ in the open argument structure above, seems to be entirely dependent (if determinate at all, as per the previous discussion of constructivity) on the possibility of associating the $\langle \mathscr{D}_i, \mathfrak{J}_i \rangle$-s to a corresponding closed $\langle \mathscr{D}, \mathfrak{J} \rangle$ for $A$ valid on $\mathfrak{C}$, under the assumption that, if those $\langle \mathscr{D}_i, \mathfrak{J}_i \rangle$-s $(i \leq n)$ exist, then such a $\langle \mathscr{D}, \mathfrak{J} \rangle$ exists too. What the $\langle \mathscr{D}, \mathfrak{J} \rangle$ is to be, is not specified by any recognisable pattern, which we could use uniformly whatever the given $\mathfrak{C}$ is. In principle, rather, there might be as many values as the $\mathfrak{C}$-s, and no schematic way to indicate all of them at once, for any ‘‘pointer". So, the reduction via $\Pi$ of the closed instances of the one-step structure that $\Pi$ is associated to, may be completely base-dependent, and there might be no rewriting rules, similar to those of $\phi_\rightarrow$ and $\iota$, capable of giving $\Pi$ an ‘‘invariant" behaviour throughout all bases.\footnote{Observe that the values may change, not only in a base-dependent way, but also in an \emph{input-dependent} way. Namely, $\Pi$ seems to be not uniform also because, \emph{even when focusing on one base}, there seems to be no schematic way to refer to the value of the reduction of a given closed instance of the one-step structure which $\Pi$ is associated to, whatever the closed instance is. This has to do, of course, with the issue touched upon in footnote 11, namely, whether uniformity should be required only relative to bases, or also relative to sets of reductions with respect to which given argument structures are valid over a given $\mathfrak{B}$. With classical logic in the meta-language, input-dependency is not a problem, as all closed $\langle \mathscr{D}_i, \mathfrak{J}_i \rangle$ for $A_i \ (i \leq n)$ valid on $\mathfrak{C} \supseteq \mathfrak{B}$, if any, can be mapped onto \emph{one and the same} closed $\langle \mathscr{D}, \mathfrak{J} \rangle$ for $A$ valid on $\mathfrak{C}$, via the constant function $\kappa_\mathfrak{C}$ of footnote 8. However, the problem of base-dependency is still there since, as soon as we put $\Pi = \bigcup_{\mathfrak{C} \supseteq \mathfrak{B}} \kappa_\mathfrak{C}$, the reductions triggered by $\Pi$ will be base-dependent, with no uniform way to refer to the values produced by each $\kappa_\mathfrak{C}$. Also, suppose we spell out in detail the constructive reading of $\forall \mathfrak{C} \supseteq \mathfrak{B} \ (\models_\mathfrak{C} \Gamma \Longrightarrow \ \models_\mathfrak{C} A)$ mentioned in footnote 8. This yields, as anticipated above, an effective function $f(x, y_1, ..., y_n)$ such that, for every $\mathfrak{C}$, $f(\mathfrak{C}, y_1, ..., y_n)$ is an effective function such that, for every closed $\langle \mathscr{D}_i, \mathfrak{J}_i \rangle$ for $A_i \ (i \leq n)$ valid on $\mathfrak{C}$ with $\Gamma = \{A_1, ..., A_n\}$, $f(\mathfrak{C}, \langle \mathscr{D}_1, \mathfrak{J}_1 \rangle, ..., \langle \mathscr{D}_n, \mathfrak{J}_n \rangle)$ is a closed $\langle \mathscr{D}, \mathfrak{J} \rangle$ for $A$ valid on $\mathfrak{C}$. Here the base-dependence is most clear, since the atomic base even occurs as a parameter of our reduction. As per Definition 18, reductions must be defined \emph{on argument structures} only. Incidentally, $f(x, y_1, ..., y_n)$ is defined on \emph{arguments}, as opposed to argument structures, i.e., on argument structures \emph{plus} sets of reductions, but this may not be a problem in itself. I will not provide the details here since I shall not be interested in this kind of solution---a solution that, after all, does not remove the dependency on the atomic base.}

Of course, non-constructivity and non-uniformity might just boil down to the same, for one could say that $\Pi$ is not constructive in the required sense just because it is not uniform in the way discussed. In fact, it is not clear whether, in the present context, constructivity and uniformity can be said to be equivalent notions, or at least whether one of them implies the other but not vice versa. In particular, while it might be sufficiently clear how the constructivity requirement should be read (for example, one could ask that reductions be recursive functions, so that the reduction sequences they generate be recursive too), it is much less clear how the uniformity requirement is to be defined more precisely.

However, providing a rigorous criterion of uniformity is beyond the scope of this paper. In Section 7, I shall single out a set of reductions which validate (logically) the atomic Split rule. All these reductions, and the set they give rise to, will turn out to be perfectly constructive and uniform, under any reasonable precisification of these notions in line with Prawitz's original intentions. For the moment, I want rather to address another important issue.

Since uniformity is not precisely defined, it might not be excluded that a result to the effect that $\Gamma \Vdash A \Longleftrightarrow \Gamma \models A$ be actually provable by choosing better strategies than $\Pi$ above for showing that $\models_\mathfrak{B}$ enjoys the same properties as $\Vdash_\mathfrak{B}$ in Definition 7. Be that as it may, it is worthwhile remarking that $\Vdash_{(\mathfrak{B})}$ differs from $\models_{(\mathfrak{B})}$ in two essential respects. These become clear if, in an attempt at establishing a link with $\models_{(\mathfrak{B})}$, we try to understand $\Gamma \Vdash_\mathfrak{B} A$, not as a primitive notion, but as meaning that there is $\langle \mathscr{D}, \mathfrak{J} \rangle$ from $\Gamma$ to $A$ which is valid on $\mathfrak{B}$.

The first and, in this context, less serious difference, has to do with logical validity. Recall that, as per Definition 8, $\Gamma \Vdash A$ is defined as $\forall \mathfrak{B} \ (\Gamma \Vdash_\mathfrak{B} A)$. If we now use our tentative interpretation above, we have that $\Gamma \Vdash A$ means

\begin{itemize}
   \item[(i)] for every $\mathfrak{B}$, there is $\langle \mathscr{D}, \mathfrak{J} \rangle$ from $\Gamma$ to $A$ which is valid over $\mathfrak{B}$.
\end{itemize}
From Definitions 23 and 24, we conclude instead that, in $\Gamma \models A$, the quantifiers of (i) are inverted, i.e., $\Gamma \models A$ means

\begin{itemize}
    \item[(ii)] there is $\langle \mathscr{D}, \mathfrak{J} \rangle$ from $\Gamma$ to $A$ such that, for every $\mathfrak{B}$, $\langle \mathscr{D}, \mathfrak{J} \rangle$ is valid on $\mathfrak{B}$.
\end{itemize}
Therefore, in Prawitz's original formulation, logical consequence does not mean that we can find a valid argument on every $\mathfrak{B}$, as obtains when interpreting $\Vdash$ in terms of existence of suitable valid arguments. Rather, it means that there is an argument which remains valid throughout all $\mathfrak{B}$-s. I will not dwell upon this issue here, as it is not problematic in a monotonic context.\footnote{Instead, \emph{it is} problematic in the non-monotonic context, since there, it boils down to basically the same issue as the one concerning what below I call the admissibility clause. See \cite{piccolominithesis, piccolominibook, piccolomininote, piccolominisomeresults} for more.} In fact, we can prove that the quantifiers in (ii) \emph{can} be inverted.

\begin{proposition}
    $\Gamma \models A \Longleftrightarrow \forall \mathfrak{B}, \exists \ \langle \mathscr{D}, \mathfrak{J} \rangle$ from $\Gamma$ to $A$ which is valid over $\mathfrak{B}$.
\end{proposition}

\begin{proof}
    ($\Longrightarrow$) is trivial. For ($\Longleftarrow$), instantiate the assumption on $\mathfrak{B}^\emptyset$. Then, there is $\langle \mathscr{D}, \mathfrak{J} \rangle$ from $\Gamma$ to $A$ valid on $\mathfrak{B}^\emptyset$. By Proposition 5, $\langle \mathscr{D}, \mathfrak{J} \rangle$ is logically valid, i.e., $\Gamma \models A$ (we could have also applied directly Corollary 2 from $\Gamma \models_{\mathfrak{B}^\emptyset} A$).
\end{proof}

The second and more serious difference, already partly mentioned above, is instead the following. Since in miB-eS $\Vdash_{(\mathfrak{B})}$ is \emph{primitive}, the condition for $\Gamma \Vdash_\mathfrak{B} A$ with $\Gamma \neq \emptyset$ to hold is \emph{defined} by the by-implication

\begin{itemize}
    \item[(**)] $\Gamma \Vdash_\mathfrak{B} A \Longleftrightarrow \forall \mathfrak{C} \supseteq \mathfrak{B} \ (\models_\mathfrak{C} \Gamma \Longrightarrow \ \models_\mathfrak{C} A)$, with $\Gamma \neq \emptyset$.
\end{itemize}
The corresponding version in miP-tV is (*) above, i.e.,

\begin{itemize}
    \item[(*)] $\Gamma \models_\mathfrak{B} A \Longleftrightarrow \ \forall \mathfrak{C} \supseteq \mathfrak{B} \ (\models_\mathfrak{C} \Gamma \Longrightarrow \ \models_\mathfrak{C} A)$, with $\Gamma \neq \emptyset$.
\end{itemize}
Since in miP-tV $\models_{(\mathfrak{B})}$ \emph{is not} primitive, but defined in terms of existence of suitable valid (over $\mathfrak{B}$) arguments from given assumptions to a given conclusion, (*) \emph{is not} part of the definition of what it is for $\Gamma \models_\mathfrak{B} A$ to hold, but \emph{must be proved} out of the definition of what it is for an open argument from $\Gamma$ to $A$ to be valid over $\mathfrak{B}$---given by Definition 21. Such a proof should establish what we obtain when translating (*) from non-primitive $\models_\mathfrak{B}$ to the primitive notion of argument valid over $\mathfrak{B}$ or, equivalently, when translating (**) via the above-told interpretation of $\Vdash_\mathfrak{B}$ in terms existence of suitable arguments valid over $\mathfrak{B}$, i.e.,

\begin{itemize}
    \item[(AC)] $\exists \ \langle \mathscr{D}, \mathfrak{J} \rangle$ valid on $\mathfrak{B}$ with $\mathscr{D}$ from $\Gamma$ to $A \Longleftrightarrow \forall \mathfrak{C} \supseteq \mathfrak{B} \ (\exists \ \langle \mathscr{D}_i, \mathfrak{J}_i \rangle$ valid on $\mathfrak{C}$ with $\mathscr{D}_i$ closed for $A_i \ (i \leq n) \Longrightarrow \exists \ \langle \mathscr{D}^*, \mathfrak{J}^* \rangle$ valid on $\mathfrak{C}$ with $\mathscr{D}^*$ closed for $A)$.
\end{itemize}
---where, as usual, $\Gamma = \{A_1, ..., A_n\}$. Now, the left-to-right direction of (AC) is, as said, easily provable in miP-tV.

\begin{proposition}
    $\Gamma \models_\mathfrak{B} A \Longrightarrow \forall \mathfrak{C} \supseteq \mathfrak{B} \ (\models_\mathfrak{C} \Gamma \Longrightarrow \ \models_\mathfrak{C} A)$.
\end{proposition}

\noindent The right-to-left direction of (AC), instead, is much bolder. It states that, if the \emph{existence} of closed arguments for each element of $\Gamma$ valid on any extension of the atomic base implies the \emph{existence} of a closed argument for $A$ valid on this extension, then there is a \emph{single} argument structure from $\Gamma$ to $A$ which be valid over the atomic base relative to some assignment of reductions. Roughly, if provability of $\Gamma$ implies provability of $A$ on any extension of a given base, then provability of $A$ \emph{from} $\Gamma$ obtains on the base---see \cite[Section 7.1.4]{piccolominibook} and \cite[Appendix]{pymcategorical} for similar observations.

It is not difficult to realise that (AC) is a sort of \emph{admissibility clause}---whence the label---and, as said, such a result is in fact provable in miP-tV when sets of reductions are allowed to be non-constructive or non-uniform. But the proof does not go through when more constraints are put on sets of reductions---a similar remark, in non-monotonic context, is raised in \cite{piccolomininote}. So, if we want to stick to the original picture of \cite{prawitz1973}, we need to find a different, possibly more direct route to incompleteness of $\texttt{IL}$.

\section{Prawitz's conjecture for miP-tV}

To stick to the original picture of \cite{prawitz1973}, I should now state what the completeness of $\texttt{IL}$ over $\models$ precisely amounts to in that picture. The issue was raised by Prawitz as a conjecture---now known as \emph{Prawitz's conjecture}---which, however, is originally not about \emph{logical consequence}, but about \emph{logically valid inference rules}.

\begin{definition}
    $\mathscr{R}$ is \emph{valid} relative to $\mathfrak{J}$ on $\mathfrak{B}$ iff, $\forall \mathfrak{H} \supseteq \mathfrak{J}$, $\forall \mathfrak{C} \supseteq \mathfrak{B}$, and for every instance $\langle \langle \mathscr{D}_1, ..., \mathscr{D}_n \rangle, A, \delta \rangle$ of $\mathscr{R}$, if $\langle \mathscr{D}_i, \mathfrak{H} \rangle \ (i \leq n)$ is valid over $\mathfrak{C}$, then
    \begin{prooftree}
        \AxiomC{$\mathscr{D}_1$}
        \AxiomC{$\dots$}
        \AxiomC{$\mathscr{D}_n$}
        \RightLabel{$\delta$}
        \TrinaryInfC{$A$}
    \end{prooftree}
    is valid over $\mathfrak{C}$ when paired with $\mathfrak{H}$.
\end{definition}

\begin{definition}
    $\mathscr{R}$ is \emph{logically valid} relative to $\mathfrak{J}$ iff $\mathscr{R}$ is valid relative to $\mathfrak{J}$ on every $\mathfrak{B}$.
\end{definition}

\begin{conjecture}[Prawitz's conjecture]
    $\mathscr{R}$ logically valid relative to some $\mathfrak{J} \Longrightarrow \ \mathscr{R}$ derivable in $\emph{\texttt{IL}}$.
\end{conjecture}
\noindent Derivability of a rule is explained by Prawitz as follows: ‘‘we may define an inference rule $\mathscr{R}$ as \emph{derivable} in a deductive system $\mathscr{S}$, when $A$ is derivable from $\Gamma$ in the extension of $\mathscr{S}$ obtained by adding the rule $\mathscr{R}$ only if $A$ is derivable from $\Gamma$ already in $\mathscr{S}$" \cite[246]{prawitz1973}. A more standard notion of completeness runs along the lines of Definition 9 above.

\begin{definition}
    $\emph{\texttt{IL}}$ is \emph{complete} over $\models$ iff $\Gamma \models A \Longrightarrow \Gamma \vdash_{\emph{\texttt{IL}}} A$.
\end{definition}

\begin{conjecture}
    $\emph{\texttt{IL}}$ is complete over $\models$.
\end{conjecture}
\noindent Before discussing the relation between Conjecture 1 and Conjecture 2, a remark is in order. This is that, for the sake of simplicity, from now on I limit myself to inference rules of the form

\begin{center}
    $\mathscr{R} = \{\langle \langle \mathscr{D}_1, ..., \mathscr{D}_n \rangle, A, \emptyset \rangle \ | \ \mathscr{D}_i$ for $A_i \ (i \leq n)\}$.
\end{center}
In other words, $\mathscr{R}$ is such that:

\begin{itemize}

\item[(a)] it can be schematically indicated by the meta-linguistic description

\begin{prooftree}
    \AxiomC{$A_1$}
    \AxiomC{$\dots$}
    \AxiomC{$A_n$}
    \TrinaryInfC{$A$}
\end{prooftree}
\item[(b)] the instances of $\mathscr{R}$ never discharge assumptions.
\end{itemize}
Both limitations serve the purpose of simplifying the proof of the equivalence between Conjecture 1 and Conjecture 2, but are otherwise harmless with respect to the incompleteness result I prove in Section 7, given that the latter concerns a rule which meets these restrictions. Limitation (a), more specifically, rules out certain limit-cases of Prawitz's notion of inference rule, e.g.: given an argument structure with conclusion $A$ and immediate sub-structures $\mathscr{D}_1, ..., \mathscr{D}_n$ which does not discharge any assumptions in the last step, the pair $\langle \langle \mathscr{D}_1, ..., \mathscr{D}_n \rangle, A, \emptyset \rangle$ is an inference by Definition 15, and the singleton $\{\langle \langle \mathscr{D}_1, ..., \mathscr{D}_n \rangle, A, \emptyset \rangle\}$ is an inference rule by Definition 16. Instead, I focus only on rules whose instances are \emph{all} the argument structures obtained by appending a certain conclusion to \emph{any} argument structures for certain premises (with no dischargement of assumptions), and which can hence be indicated by schematic meta-linguistics descriptions where only those premises and that conclusion occur (without indications of dischargement of assumptions).

\begin{proposition}
    With $\Gamma = \{A_1, ..., A_n\}$, $\Gamma \models A$ iff the inference rule
    \begin{center}
        $\mathscr{R} = \{\langle \langle \mathscr{D}_1, ..., \mathscr{D}_n \rangle, A, \emptyset \rangle \ | \ \mathscr{D}_i$ for $A_i \ (i \leq n)\}$
    \end{center}
    is logically valid relative to some $\mathfrak{J}$.
\end{proposition}

\begin{proof}
    ($\Longrightarrow$) Assume $\Gamma \models A$, i.e., there is $\langle \mathscr{D}, \mathfrak{J} \rangle$ logically valid with $\mathscr{D}$ from $\Gamma$ to $A$, written
    \begin{prooftree}
        \AxiomC{$A_1$}
        \AxiomC{$\dots$}
        \AxiomC{$A_n$}
        \noLine
        \TrinaryInfC{$\mathscr{D}$}
        \noLine
        \UnaryInfC{$A$}
    \end{prooftree}
\noindent Consider the function $\phi$ defined on $\mathscr{R}$ such that, for every $\mathscr{D}_i$ for $A_i \ (i \leq n)$,
    \begin{prooftree}
        \AxiomC{$\mathscr{D}_1$}
        \noLine
        \UnaryInfC{$A_1$}
        \AxiomC{$\dots$}
        \AxiomC{$\mathscr{D}_n$}
        \noLine
        \UnaryInfC{$A_n$}
        \TrinaryInfC{$A$}
        \AxiomC{$\stackrel{\phi}{\Longrightarrow}$}
        \noLine
        \UnaryInfC{}
        \AxiomC{$\mathscr{D}_1$}
        \noLine
        \UnaryInfC{$A_1$}
        \AxiomC{$\dots$}
        \AxiomC{$\mathscr{D}_n$}
        \noLine
        \UnaryInfC{$A_n$}
        \noLine
        \TrinaryInfC{$\mathscr{D}$}
        \noLine
        \UnaryInfC{$A$}
        \noLine
        \TrinaryInfC{}
    \end{prooftree}
    It is easy to see that $\phi$ is constructive and uniform, and that it satisfies points 1 to 4 in Definition 18, so $\phi$ is a reduction for $\mathscr{R}$ (in the strict understanding of miP-tV). Now, $\mathscr{R}$ is logically valid relative to $\mathfrak{J} \cup \{\phi\}$. Indeed, suppose that $\langle \mathscr{D}_i, \mathfrak{H} \rangle$ with $\mathscr{D}_i$ for $A_i \ (i \leq n)$ and $\mathfrak{H} \supseteq \mathfrak{J} \cup \{\phi\}$ is valid over any $\mathfrak{C} \supseteq \mathfrak{B}$ of any $\mathfrak{B}$. We can safely assume that the $\mathscr{D}_i$-s $(i \leq n)$ are closed. We have to show that $\langle \mathscr{D}^*, \mathfrak{H} \rangle$ is also valid over $\mathfrak{C}$, where $\mathscr{D}^*$ is
    \begin{prooftree}
        \AxiomC{$\mathscr{D}_1$}
        \noLine
        \UnaryInfC{$A_1$}
        \AxiomC{$\dots$}
        \AxiomC{$\mathscr{D}_n$}
        \noLine
        \UnaryInfC{$A_n$}
        \TrinaryInfC{$A$}
    \end{prooftree}
    Since $\phi \in \mathfrak{H}$, we have $\mathscr{D}^* \leq_\mathfrak{H} \mathscr{D}^{**}$, where $\mathscr{D}^{**}$ is
    \begin{prooftree}
        \AxiomC{$\mathscr{D}_1$}
        \noLine
        \UnaryInfC{$A_1$}
        \AxiomC{$\dots$}
        \AxiomC{$\mathscr{D}_n$}
        \noLine
        \UnaryInfC{$A_n$}
        \noLine
        \TrinaryInfC{$\mathscr{D}$}
        \noLine
        \UnaryInfC{$A$}
    \end{prooftree}
    Since $\langle \mathscr{D}, \mathfrak{J} \rangle$ is logically valid, by Definitions 21 and 23, $\langle \mathscr{D}^{**}, \mathfrak{H} \rangle$ is valid over $\mathfrak{C}$, hence $\mathscr{D}^{**} \leq_\mathfrak{H} \mathscr{D}^{***}$, where $\mathscr{D}^{***}$ is a closed canonical argument structure for $A$, whose immediate sub-structures are valid on $\mathfrak{C}$ relative to $\mathfrak{H}$. By definition 20, $\leq_\mathfrak{H}$ is transitive, so we have $\mathscr{D}^* \leq_\mathfrak{H} \mathscr{D}^{***}$. ($\Longleftarrow$) Suppose $\mathscr{R}$ is logically valid relative to some $\mathfrak{J}$, and consider the one-step argument structure from $A_1, ..., A_n$ to $A$, which we may indicate by $\mathscr{D}$. We have that $\langle \mathscr{D}, \mathfrak{J} \rangle$ is logically valid. For, take any $\langle \mathscr{D}_i, \mathfrak{H} \rangle$ valid over any $\mathfrak{C} \supseteq \mathfrak{B}$ of any $\mathfrak{B}$, with $\mathscr{D}_i$ closed for $A_i \ (i \leq n)$ and $\mathfrak{H} \supseteq \mathfrak{J}$. The argument structure
    \begin{prooftree}
        \AxiomC{$\mathscr{D}_1$}
        \noLine
        \UnaryInfC{$A_1$}
        \AxiomC{$\dots$}
        \AxiomC{$\mathscr{D}_n$}
        \noLine
        \UnaryInfC{$A_n$}
        \TrinaryInfC{$A$}
    \end{prooftree}
    which we may indicate by $\mathscr{D}^*$, is such that $\mathscr{D}^* \in \mathscr{R}$. Hence, by Definition 25, $\langle \mathscr{D}^*, \mathfrak{H} \rangle$ is valid over $\mathfrak{C}$.
\end{proof}

\begin{proposition}
    Conjecture 1 holds iff Conjecture 2 holds.
\end{proposition}

\begin{proof}
    ($\Longrightarrow$) Assume $\Gamma \models A$, with $\Gamma = \{A_1, ..., A_n\}$. Then, by Proposition 10, the inference rule
    \begin{center}
        $\mathscr{R} = \{\langle \langle \mathscr{D}_1, ..., \mathscr{D}_n \rangle, A, \emptyset \rangle \ | \ \mathscr{D}_i$ for $A_i \ (i \leq n)\}$
    \end{center}
    is logically valid with respect to some $\mathfrak{J}$. But since we assumed Conjecture 1 to hold, $\mathscr{R}$ is derivable in $\texttt{IL}$, i.e., $\Gamma \vdash_\texttt{IL} A$. As for ($\Longleftarrow$), recall that I am restricting myself to inference rules with limitations (a) and (b). Assume that such a rule

    \begin{center}
        $\mathscr{R} = \{\langle \langle \mathscr{D}_1, ..., \mathscr{D}_n \rangle, A, \emptyset \rangle \ | \ \mathscr{D}_i$ for $A_i \ (i \leq n)\}$.
    \end{center}
    with $\Gamma = \{A_1, ..., A_n\}$ is logically valid. By Proposition 10, $\Gamma \models A$ and, by Conjecture 2, $\Gamma \vdash_\texttt{IL} A$.
\end{proof}

\noindent It is important to remark that, if we decide to formulate the completeness conjecture in the terms of Conjecture 2, we should not forget to demand---if we aim at sticking to the original picture of \cite{prawitz1973}---that the logically valid argument required for $\Gamma \models A$ to hold involves a \emph{constructive} set of \emph{uniform} reductions, as per Section 4. As said, Conjecture 2 \emph{is not} Prawitz's original conjecture, which was given in the form of Conjecture 1 instead. And in the latter it is understood that logical validity of inferences rules should be relative to constructive sets of uniform reductions.

\section{Incompleteness via constructions}

Especially interesting, for the purposes of this paper, are a number of results proved (with essentially one and the same strategy) in \cite{piechadecampossanz, piechadecampossanzschroederheisterpeirce}. These establish the incompleteness of $\texttt{IL}$ relative to a notion of logical consequence drawn from an approach developed by Prawitz in \cite{prawitz1971}. The latter is in turn based on a concept of construction inspired by the concept of proof involved in BHK semantics \cite{troelstravandalen}.

Assuming the underlying notions of atomic rule, atomic base, and atomic derivation given in Definitions  2, 3 and 5 above, let us state what follows---I limit myself to the $\wedge$-free fragment of $\mathscr{L}$, which is enough for my purposes here.

\begin{definition}
    $k$ is a \emph{construction of} $A \in \emph{\texttt{FORM}}_\mathscr{L}$ \emph{on} $\mathfrak{B}$ iff
    \begin{itemize}
        \item $A \in \emph{\texttt{ATOM}}_\mathscr{L} \Longrightarrow \ k \in \emph{\texttt{DER}}_\mathfrak{B}$ witnesses $\vdash_\mathfrak{B} A$
        \item $A = B_1 \vee B_2 \Longrightarrow \ k$ is $\langle i, k^* \rangle$, with $i = 1, 2$ and $k^*$ construction of $B_i$ on $\mathfrak{B}$
        \item $A = B \rightarrow C \Longrightarrow \ k$ is $\lambda B.k^*(B)$, with $k^*(B)$ construction of $C$ from $B$ on $\mathfrak{B}$.
    \end{itemize}
\end{definition}

\begin{definition}
    $k(A_1, ..., A_n)$ is a \emph{construction of} $A$ \emph{from} $\{A_1, ..., A_n\}$ on $\mathfrak{B}$ iff, for every $\mathfrak{C} \supseteq \mathfrak{B}$, if $k_i$ is a construction of $A_i \ (i \leq n)$ on $\mathfrak{C}$, then $k(k_1, ..., k_n)$ is a construction of $A$ on $\mathfrak{C}$.
\end{definition}

\begin{proposition}
    Let $\mathscr{D}$ be from $p_1, ..., p_n$ to $q$, and suppose that, when assumptions $p_1, ..., p_n$ are turned into axioms, $\mathscr{D}$ becomes an atomic derivation in the extension of $\mathfrak{B}$ obtained by adding these axioms to it. Then, $\mathscr{D}$ is a construction of $q$ from $p_1, ..., p_n$ on $\mathfrak{B}$.
\end{proposition}

\begin{proof}
    Consider any $\mathfrak{C} \supseteq \mathfrak{B}$, and take any construction $k_i$ of $p_i \ (i \leq n)$ on $\mathfrak{C}$. Thus, $k_i \in \texttt{DER}_\mathfrak{C}$ witnesses $\vdash_\mathfrak{C} p_i \ (i \leq n)$. But then, $\mathscr{D}[k_1, ..., k_n/p_1, ..., p_n] \in \texttt{DER}_\mathfrak{C}$ witnesses $\vdash_\mathfrak{C} q$, so it is a construction of $q$ on $\mathfrak{C}$.
\end{proof}

\begin{definition}
    $A$ is a \emph{logical consequence of} $\Gamma$ \emph{in Prawitz's 1971 sense}---written $\Gamma \models^c A$---iff $\exists k, \forall \mathfrak{B}$, $k$ is a construction of $A$ from $\Gamma$ on $\mathfrak{B}$.
\end{definition}

\noindent The following proof for the logical validity of the atomic Split rule relative to $\models^c$ is from \cite{piechadecampossanz}, with the exception that, there, the result concerns Mints rule---a conceptually similar proof is \cite[Theorem 37]{pymcategorical}. Logical validity of the atomic Split rule is proved in \cite{goldfarb, piechaschroeder-heisterdecampossanz} but, although the proof-strategy is essentially the same as in \cite{piechadecampossanz} and \cite{piechadecampossanzschroederheisterpeirce}, the proof has an accent which is less in line than \cite{piechadecampossanz} relative to my aims here. The proof in \cite{piechadecampossanzschroederheisterpeirce} is for the logical validity of Peirce rule, but---as the authors observe---it works only when atomic bases have level at most $1$, so it is not suitable for the kind of approach I am pursuing here.

\begin{theorem}
    $p \rightarrow q \vee r \models^c (p \rightarrow q) \vee (p \rightarrow r)$.
\end{theorem}

\begin{proof}
    We have to show that there is $k$ such that, for every $\mathfrak{B}$, for every $\mathfrak{C} \supseteq \mathfrak{B}$, if $k^*$ is a construction of $p \rightarrow q \vee r$ on $\mathfrak{C}$, then $k(k^*)$ is a construction of $(p \rightarrow q) \vee (p \rightarrow r)$ on $\mathfrak{C}$. Take any $\mathfrak{B}$ and any $\mathfrak{C} \supseteq \mathfrak{B}$, and suppose $\lambda p. k_1(p)
    $ is a construction of $p \rightarrow q \vee r$ on $\mathfrak{C}$. So, for any $\mathfrak{D} \supseteq \mathfrak{C}$, if $k^{**}$ is a construction of $p$ on $\mathfrak{D}$, $k_1(k^{**})$ is a construction of $q \vee r$ on $\mathfrak{D}$. Consider now the extension $\mathfrak{D}$ of $\mathfrak{C}$ given by adding to $\mathfrak{C}$ the axiom

    \begin{prooftree}
        \AxiomC{}
        \RightLabel{$\mathfrak{p}$}
        \UnaryInfC{$p$}
    \end{prooftree}
    An instance of $\mathfrak{p}$ yields $k_2 \in \texttt{DER}_\mathfrak{D}$ witnessing $\vdash_\mathfrak{D} p$, and $k_1(k_2)$ is a construction of $q \vee r$ in $\mathfrak{D}$. So: either $k_1(k_2) = \langle 1, k_3 \rangle$, with $k_3$ construction of $q$ on $\mathfrak{D}$; or $k_1(k_2) = \langle 2, k_3 \rangle$, with $k_3$ construction of $r$ on $\mathfrak{D}$. In both cases we have $k_3 \in \texttt{DER}_\mathfrak{D}$; in the first case, $k_3$ witnesses $\vdash_\mathfrak{D} q$ whereas, in the second case, $k_3$ witnesses $\vdash_\mathfrak{D} r$. If no instances of $\mathfrak{p}$ are in $k_3$, $k_3 \in \texttt{DER}_\mathfrak{C}$. But then $\lambda p.k_3$ is (vacuously) a construction of $q$ [resp. $r$] from $p$, whence $\langle 1, \lambda p.k_3 \rangle$ [resp. $\langle 2, \lambda p.k_3 \rangle$] is a construction of $(p \rightarrow q) \vee (p \rightarrow r)$ on $\mathfrak{C}$. If instead $k_3$ contains some instances of $\mathfrak{p}$, we have that $k_3[p/\mathfrak{p}]$---i.e., the result of replacing the application of any instance of $\mathfrak{p}$ in $k_3$ with $p$ as an assumption---is an argument structure from $p$ to $q$ [resp. from $p$ to $r$] which, by Proposition 12, is a construction of $q$ [resp. $r$] from $p$ on $\mathfrak{C}$. So, $\langle 1, \lambda p.k_3[p/\mathfrak{p}] \rangle$ [resp. $\langle 2, \lambda p.k_3[p/\mathfrak{p}]] \rangle$] is a construction of $(p \rightarrow q) \vee (p \rightarrow r)$ on $\mathfrak{C}$. The procedure of extending $\mathfrak{C}$ via axiom $\mathfrak{p}$, and then looking for a derivation of $(p \rightarrow q) \vee (p \rightarrow r)$ is the desired construction $k$.
\end{proof}

Contrarily to an attempted translation from miB-eS to miP-tV, Theorem 2 provides a concrete witness of the logical validity of the atomic Split rule. So, based on what I said in Section 5, this result, and the way it is proved, are more in line with the spirit of \cite{prawitz1973}. On the other hand, however, Theorem 2 concerns $\models^c$, namely, a notion of proof-theoretic (logical) consequence based on the notion of (intuitionistic) construction. Thus, in asking whether something similar is available for $\models$ too, we should investigate whether the construction $k$ in the proof of Theorem 2 can be turned into a reduction or set of reductions which be constructive and uniform in the sense discussed in Section 4.

One may wonder about the precise relation between (intuitionistic) constructions, and the kind of reductions or sets of reductions at play in miP-tV. The idea seems to be that the latter are instances of the former, while the inverse fails:

\begin{quote}
    an implication $A \rightarrow B$ is constructively understood as the assertion of the existence of a construction of $B$ from $A$, and in accordance with this meaning $\rightarrow_I$ allows the inference of $A \rightarrow B$ given a proof of $B$ from $A$. Of course, a proof of $B$ from $A$ is not the same as a construction of $B$ from $A$; it is rather a special kind of such a construction. \cite[244]{prawitz1971}
\end{quote}
Prawitz also says that a proof of $B$ from $A$ ‘‘gives a particular kind of a uniform construction transforming constructions of $A$ to constructions of $B$, which is not required by the usual constructive interpretation" \cite[276]{prawitz1971}. Schroeder-Heister has instead remarked that sets of reductions in Prawitz's sense undergo what he calls a \emph{full specification} requirement:

\begin{quote}
    we should require of a justification $\mathfrak{J}$ that it contains \emph{all} reductions, which any reduction in $\mathfrak{J}$ refers to. [...]. I consider this to be crucial for the proof-theoretic approach in contradistinction to an approach which may be called \emph{constructive}. \cite[Appendix, 4]{schroeder-heisterbasicideas}
\end{quote}
Therefore, sets of reductions are partial constructive functions of a special kind, the main difference being that ‘‘for an arbitrary partial constructive functional $f$ we would allow that as a value it produces a partial function(al) $f'$ which is defined on some arguments for which $f$ is not defined" \cite[Appendix, 8]{schroeder-heisterbasicideas}.

Let us then go back to the construction $k$ in the proof of Theorem 2, and let us indicate by $\pi_1$ [resp. $\pi_2$] the standard left [resp. right] projection on a pair, and by $\mathfrak{p} \in \mathscr{D}$ the fact that $\mathfrak{p}$ is used in $\mathscr{D}$, where $\mathfrak{p}$ is as in the proof of Theorem 2. Our $k$ can be defined by cases as follows:

    \begin{center}
        $k(k_1) = \begin{cases} \langle \pi_1(k_1(k_2)), \lambda p.\pi_2(k_1(k_2)) \rangle & \mathfrak{p} \notin \pi_2(k_1(k_2)) \\ \langle \pi_1(k_1(k_2)), \lambda p.((\pi_2(k_1(k_2)))[p/\mathfrak{p}]) \rangle & \text{otherwise} \end{cases}$
    \end{center}
where $k_1$ is a construction for $p \rightarrow q \vee r$ as in the proof of Theorem 2. Now, \emph{qua} construction, $k$ seems to be perfectly uniform, as it behaves in the same way over all bases. But since, for what said above, not all intuitionistic constructions are sets of reductions as required by miP-tV, it is not clear whether from $k$  we can extract a set $\mathfrak{J}$ of uniform reductions such that the atomic Split rule is logically valid over $\models$.

For example, we may not exclude \emph{a priori} that any set of reductions for the atomic Split rule is non-uniform. This is, e.g., what happens with immediate attempts at turning $k$ into a Prawitzian reduction. The latter should operate starting from argument structures

\begin{prooftree}
    \AxiomC{$[p]_1$}
    \noLine
    \UnaryInfC{$\mathscr{D}_1$}
    \noLine
    \UnaryInfC{$q \vee r$}
    \RightLabel{$1$}
    \UnaryInfC{$p \rightarrow q \vee r$}
    \LeftLabel{$\mathscr{D} = \ $}
    \UnaryInfC{$(p \rightarrow q) \vee (p \rightarrow r)$}
\end{prooftree}
Let us assume, paralleling $k_1$ in the proof of Theorem 2, that $\mathscr{D}_1$ above is valid on $\mathfrak{C} \supseteq \mathfrak{B}$ relative to a set of reductions $\mathfrak{J}$. When assumption $p$ is replaced in $\mathscr{D}_1$ by $p$ used as an axiom, the resulting structure reduces to a closed atomic derivation $\mathscr{D}_2 \in \texttt{DER}_\mathfrak{D}$ for either $q$ or $r$ where, as in the proof of Theorem 2, $\mathfrak{D} = \mathfrak{C} \cup \{\mathfrak{p}\}$. If $\mathscr{D}_2 \in \texttt{DER}_\mathfrak{C}$, we  apply $\rightarrow_I$ vacuously on it, and then we apply introduction of $\vee$, getting a closed argument structure for $(p \rightarrow q) \vee (p \rightarrow r)$ which is valid on $\mathfrak{C}$ relative to $\mathfrak{J}$. If instead $\mathscr{D}_2 \notin \texttt{DER}_\mathfrak{C}$, we replace $p$ as an axiom by $p$ as an assumption, and get an open $\mathscr{D}_5$ which (as proved below) is valid over $\mathfrak{C}$ relative to $\mathfrak{J}$. We use $\rightarrow_I$ again---not vacuously but discharging assumption $p$---and introduction of $\vee$, and get again a closed argument structure for $(p \rightarrow q) \vee (p \rightarrow r)$ valid over $\mathfrak{C}$ relative to $\mathfrak{J}$.

So, whenever we are given a closed $\langle \mathscr{D}_1, \mathfrak{J} \rangle$ for $p \rightarrow q \vee r$ valid over any extension $\mathfrak{C}$ of a given atomic base $\mathfrak{B}$, we can find a closed $\langle \mathscr{D}^*, \mathfrak{J} \rangle$ for $(p \rightarrow q) \vee (p \rightarrow r)$ valid over $\mathfrak{C}$. We may thus think of a reduction

\begin{prooftree}
    \AxiomC{$[p]_1$}
    \noLine
    \UnaryInfC{$\mathscr{D}_1$}
    \noLine
    \UnaryInfC{$q \vee r$}
    \RightLabel{$1$}
    \UnaryInfC{$p \rightarrow q \vee r$}
    \UnaryInfC{$(p \rightarrow q) \vee (p \rightarrow r)$}
    \AxiomC{$\stackrel{\phi_{\texttt{Split}}}{\Longrightarrow}$}
    \noLine
    \UnaryInfC{}
    \AxiomC{$\mathscr{D}^*$}
    \noLine
    \UnaryInfC{}
    \noLine
    \TrinaryInfC{}
\end{prooftree}
where, according to the procedure above, $\mathscr{D}^*$ is described as the structure obtained by applying $\rightarrow_I$ and introduction of $\vee$ to the structure that $\mathscr{D}_1$ reduces to when closed with axiom $p$ relative to the extension obtained by adding this axiom to the extension of the $\mathfrak{B}$ which $\mathscr{D}_1$ is valid on---after possibly turning $p$ into assumption again.

It is clear, though, that $\phi_{\texttt{Split}}$ is not in line with the original intentions of \cite{prawitz1973}, as described in Section 4. In particular, it is non-uniform, as the value $\mathscr{D}^*$ cannot be expected to be the same, or to have an invariant schematic form, over all (extensions of) atomic bases over which $p \rightarrow q \vee r$ is valid.\footnote{$\phi_{\texttt{Split}}$ is non-uniform also when uniformity is required for sets of reductions, because what $\mathscr{D}^*$ is depends, not only on (the extension of) the given base, but also on the set of reductions relative to which $\mathscr{D}_1$ is valid on that (extension of the) base. $\phi_{\texttt{Split}}$ is not deterministic either, as one and the same $\mathscr{D}_1$ might be valid, relative to the same $\mathfrak{J}$, on two distinct (extensions of) atomic bases, but reduce to different values on each.}

\section{Uniform reductions for atomic Split rule}

The base dependencies of $\phi_{\texttt{Split}}$ can as said be removed. Building upon the result discussed in Section 6, one can in fact identify a set of reductions $\mathfrak{J}$ which, as required by \cite{prawitz1973}, is constructive and uniform under any precisification of the notions of constructivity and uniformity, and such that the rule

\begin{prooftree}
    \AxiomC{$p \rightarrow q \vee r$}
    \RightLabel{$\texttt{Split}$}
    \UnaryInfC{$(p \rightarrow q) \vee (p \rightarrow r)$}
\end{prooftree}
is logically valid relative to $\mathfrak{J}$. Before turning to this, let me first state a result similar to Proposition 12.

\begin{proposition}
  Let $\mathscr{D}$ be from $p_1, ..., p_n$ to $q$, and suppose that, when assumptions $p_1, ..., p_n$ are turned into axioms, $\mathscr{D}$ becomes an atomic derivation in the extension of $\mathfrak{B}$ obtained by adding these axioms to it. Then, $\langle \mathscr{D}, \emptyset \rangle$ is valid on $\mathfrak{B}$.   
\end{proposition}

\begin{proof}
    Our $\mathscr{D}$ has the form

    \begin{prooftree}
        \AxiomC{$p_1$}
        \AxiomC{$\dots$}
        \AxiomC{$p_n$}
        \noLine
        \TrinaryInfC{$\mathscr{D}$}
        \noLine
        \UnaryInfC{$q$}
    \end{prooftree}
    We have to show that, for every $\mathfrak{C} \supseteq \mathfrak{B}$, every $\mathfrak{J}$, and every closed $\langle \mathscr{D}_i, \mathfrak{J} \rangle$ for $p_i \ (i \leq n)$ valid on $\mathfrak{C}$,
    \begin{prooftree}
        \AxiomC{$\mathscr{D}_1$}
        \noLine
        \UnaryInfC{$p_1$}
        \AxiomC{$\dots$}
        \AxiomC{$\mathscr{D}_n$}
        \noLine
        \UnaryInfC{$p_n$}
        \noLine
        \TrinaryInfC{$\mathscr{D}$}
        \noLine
        \UnaryInfC{$q$}
    \end{prooftree}
    is also valid on $\mathfrak{C}$ relative to $\mathfrak{J}$. By assumption, $\mathscr{D}_i \leq_\mathfrak{J} \mathscr{D}^*_i$ with $\mathscr{D}^*_i \in \texttt{DER}_\mathfrak{C}$. But then we have
    \begin{center}
        $\mathscr{D}[\mathscr{D}_1, ..., \mathscr{D}_n/p_1, ..., p_n] \leq_\mathfrak{J} \mathscr{D}[\mathscr{D}^*_1, ..., \mathscr{D}_n/p_1, ..., p_n] \leq_\mathfrak{J} ... \leq_\mathfrak{J} \mathscr{D}[\mathscr{D}^*_1, ..., \mathscr{D}^*_n/p_1, ..., p_n]$
    \end{center}
    where $\mathscr{D}[\mathscr{D}^*_1, ..., \mathscr{D}^*_n/p_1, ..., p_n] \in \texttt{DER}_\mathfrak{C}$.
\end{proof}

\begin{corollary}
    Let $\mathscr{D}$ be as in Proposition 13. Then $\langle \mathscr{D}, \mathfrak{J} \rangle$ is valid on $\mathfrak{B}$ for every $\mathfrak{J}$.
\end{corollary}

A key-observation is now in order towards removing the base-dependencies of $\phi_\texttt{Split}$ from Section 6. To begin with, observe that argument structures are mere syntactic objects, given independently of the criteria by which they are valid relative to some $\mathfrak{J}$ on some $\mathfrak{B}$ (Definition 10). Likewise, reductions are just mappings from and to argument structures, whose definition, hence outputs they generate from given inputs, are also independent of the potential validity of those inputs or outputs relative to some $\mathfrak{J}$ over some $\mathfrak{B}$ (Definition 18). It follows that, given some $\mathfrak{J}$, the reduction sequences modulo $\mathfrak{J}$---i.e., the sequences $\mathscr{D}_1, ..., \mathscr{D}_n$ such that $\mathscr{D}_i \leq_\mathfrak{J} \mathscr{D}_{i + 1}$---are available on \emph{all} $\mathfrak{B}$-s, independently of whether the $\mathscr{D}_i$-s $(i \leq n)$ are or not valid on such $\mathfrak{B}$-s relative to $\mathfrak{J}$, or to other sets of reductions. Hence, the question whether $\mathscr{D}$ reduces to $\mathscr{D}^*$ modulo $\mathfrak{J}$, and the question whether $\mathscr{D}$ or $\mathscr{D}^*$ are valid relative to $\mathfrak{J}$ on some $\mathfrak{B}$, are kept apart from each other. They may of course be connected, as the validity of $\mathscr{D}$ relative to $\mathfrak{J}$ on $\mathfrak{B}$ might happen to stem from the reducibility of $\mathscr{D}$ to $\mathscr{D}^*$ modulo $\mathfrak{J}$, with $\mathscr{D}^*$ valid on $\mathfrak{B}$ relative to $\mathfrak{J}$. When this obtains, however, we are in the presence of a conjunction whose conjuncts concern different phenomena, i.e., $\mathscr{D} \leq_\mathfrak{J} \mathscr{D}^*$, on the one hand, and $\langle \mathscr{D}^*, \mathfrak{J} \rangle$ valid on $\mathfrak{B}$ on the other.
Hence, to prove the validity of an argument structure relative to $\mathfrak{J}$ on $\mathfrak{B}$, we can make ‘‘detours" through any $\mathfrak{C}$ different from $\mathfrak{B}$, relying only on the reduction sequences triggered by $\mathfrak{J}$, and without caring about whether argument structures involved in these sequences are or not valid on $\mathfrak{B}$ relative to $\mathfrak{J}$. The reduction sequences will be in any case available, just because of how the reductions at issue are defined.

Armed with this remark, which is used implicitly in the proof of Theorem 3 below, let us now go back to $\texttt{Split}$, and let us define two mappings. First, a mapping $\phi^1$ defined on the set

\begin{center}
    $\mathbb{D}^1 = \{\langle \mathscr{D}, (p \rightarrow q) \vee (p \rightarrow r), \emptyset \rangle \ | \ \begin{rcases} \mathscr{D} \in \ \rightarrow_I \ \text{closed for} \ p \rightarrow q \vee r\\\mathscr{D} \ \text{has immediate sub-argument} \ \mathscr{D}^* \ \text{open}\end{rcases}$
\end{center}
which operates as follows

\begin{prooftree}
    \AxiomC{$[p]_1$}
    \noLine
    \UnaryInfC{$\mathscr{D}^*$}
    \noLine
    \UnaryInfC{$q \vee r$}
    \RightLabel{$1$}
    \UnaryInfC{$p \rightarrow q \vee r$}
    \RightLabel{$\texttt{Split}$}
    \UnaryInfC{$(p \rightarrow q) \vee (p \rightarrow r)$}
    \AxiomC{$\stackrel{\phi^1}{\Longrightarrow}$}
    \noLine
    \UnaryInfC{}
    \AxiomC{}
    \UnaryInfC{$p$}
    \noLine
    \UnaryInfC{$\mathscr{D}^{**}$}
    \noLine
    \UnaryInfC{$q \vee r$}
    \UnaryInfC{$p \rightarrow q \vee r$}
    \RightLabel{$\texttt{Split}$}
    \UnaryInfC{$(p \rightarrow q) \vee (p \rightarrow r)$}
    \noLine
    \TrinaryInfC{}
\end{prooftree}
Thus, $\phi^1$ turns a closed canonical argument structure for $p \rightarrow q \vee r$ whose immediate sub-argument is open into one where assumption $p$ is turned into an axiom, whilst $p$ is discharged vacuously in the step of $\rightarrow_I$ yielding the premise of the last inference. The rest stays unchanged---i.e., $\mathscr{D}^{**}$ is in all ways identical to $\mathscr{D}^*$, except for the changes indicated. Next we have a mapping $\phi^2$ defined on the set

\begin{center}
    $\mathbb{D}^2 = \{\langle \mathscr{D}, (p \rightarrow q_1) \vee (p \rightarrow q_2), \emptyset \rangle \ | \ \begin{rcases} \mathscr{D} \in \ \rightarrow_I \ \text{closed for} \ p \rightarrow q_1 \vee q_2\\\mathscr{D} \ \text{has immediate sub-argument} \ \mathscr{D}^* \ \text{canonical}\\ \text{the immediate sub-argument} \ \mathscr{D}^{**} \ \text{of} \ \mathscr{D}^* \ \text{is an atomic derivation}\end{rcases}$
\end{center}
which operates as follows

\begin{prooftree}
    \AxiomC{}
    \UnaryInfC{$p$}
    \noLine
    \UnaryInfC{$\mathscr{D}^{**}$}
    \noLine
    \UnaryInfC{$q_i$}
    \UnaryInfC{$q_1 \vee q_2$}
    \UnaryInfC{$p \rightarrow q_1 \vee q_2$}
    \RightLabel{$\texttt{Split}$}
    \UnaryInfC{$(p \rightarrow q_1) \vee (p \rightarrow q_2)$}
    \AxiomC{$\stackrel{\phi^2}{\Longrightarrow}$}
    \noLine
    \UnaryInfC{}
    \AxiomC{$[p]_1$}
    \noLine
    \UnaryInfC{$\mathscr{D}^{***}$}
    \noLine
    \UnaryInfC{$q_i$}
    \RightLabel{$1$}
    \UnaryInfC{$p \rightarrow q_i$}
    \UnaryInfC{$(p \rightarrow q_1) \vee (p \rightarrow q_2)$}
    \noLine
    \TrinaryInfC{}
\end{prooftree}
Hence, $\phi^2$ turns a closed canonical argument structure for $p \rightarrow q_1 \vee q_2$ where $p$ is vacuously discharged by the step of $\rightarrow_I$ yielding the premise of the last inference, and whose immediate sub-structure is in turn canonical with immediate sub-structure atomic where axiom $p$ may occur (undischarged), into one where, first, any undischarged application of axiom $p$ is turned into an assumption and, second, the application of introduction of disjunction is swapped with that of $\rightarrow_I$, while discharging assumption $p$ (not vacuously if axiom $p$ was used undischarged in the original argument structure).

Observe that $\mathbb{D}^1 \cap \mathbb{D}^2 = \emptyset$. For, suppose $\mathscr{D} \in \mathbb{D}^1$. Now, the immediate sub-structure $\mathscr{D}^*$ of $\mathscr{D}$ is canonical, and the immediate sub-structure $\mathscr{D}^{**}$ of $\mathscr{D}^*$ is open, so the only undischarged assumption of $\mathscr{D}^{**}$ must be $p$. So, even assuming that $\mathscr{D}^{**}$ is canonical, its immediate sub-structure cannot be an atomic derivation since, by Definition 5, in the latter we can have no undischarged assumptions (but only undischarged rules). Hence, $\phi^1$ and $\phi^2$ are deterministic.

\begin{proposition}
    $\phi^i$ is a reduction for $\mathbb{D}^i \ (i = 1, 2)$ according to Definition 18. 
\end{proposition}

\begin{proof}
    We have to check that $\phi^i$ satisfies the conditions 1, 2, 3 and 4 from Definition 18.
    \begin{enumerate}
        \item It is easy to check that $\phi^i$ goes from argument structures to argument structures.
        \item It is also easy to check that $\phi^i$ is defined on a subset $\mathbb{D} \subseteq \mathbb{D}^i$, i.e., $\mathbb{D} = \mathbb{D}^i$. Also, since every $\mathscr{D} \in \mathbb{D}^i$ is closed, we have $\mathscr{D}^\sigma = \mathscr{D}$ for every $\sigma$, hence trivially, $\forall \sigma, \mathscr{D} \in \mathbb{D}^i \Longrightarrow \mathscr{D}^\sigma \in \mathbb{D}^i$.
        \item Since every $\mathscr{D} \in \mathbb{D}^i$ is closed, we have just to check that, for every $\mathscr{D} \in \mathbb{D}^i$, $\phi^i(\mathscr{D})$ is closed---since clearly $\phi^i(\mathscr{D})$ has the same conclusion as $\mathscr{D}$. Now, $\phi^i(\mathscr{D})$ is obtained from $\mathscr{D}$ without adding further assumptions, and by either turning a discharged assumption into an axiom, or vice versa, or else by keeping discharged all the assumptions which were discharged in $\mathscr{D}$. Whence, $\phi^i(\mathscr{D})$ is closed.
        \item Since every $\mathscr{D} \in \mathbb{D}^i$ is closed, and since $\phi^i(\mathscr{D})$ is closed for every $\mathscr{D} \in \mathbb{D}^i$, we have $\mathscr{D}^\sigma = \mathscr{D}$ and $\phi(\mathscr{D})^\sigma = \phi(\mathscr{D})$ for every $\sigma$, whence $\phi(\mathscr{D}^\sigma) = \phi(\mathscr{D}) = \phi(\mathscr{D})^\sigma$ for every $\sigma$, i.e., $\forall \sigma, \phi(\mathscr{D}^\sigma) = \phi(\mathscr{D})^\sigma$.
    \end{enumerate}
    It is also easy to verify that the reductions are constructive and uniform as (implicitly) required in miP-tV.
\end{proof}

\begin{theorem}
    $\emph{\texttt{Split}}$ is logivally valid relative to $\mathfrak{H} = \{\phi^1, \phi^2\}$.
\end{theorem}

\begin{proof}
    We must prove that, given any $\mathfrak{B}$, any $\mathfrak{C} \supseteq \mathfrak{B}$, and any $\langle \mathscr{D}_1, \mathfrak{J} \rangle$ for $p \rightarrow q_1 \vee q_2$ valid on $\mathfrak{C}$, the argument structure

    \begin{prooftree}
        \AxiomC{$\mathscr{D}_1$}
        \noLine
        \UnaryInfC{$p \rightarrow q_1 \vee q_2$}
        \LeftLabel{$\mathscr{D}_2 = $}
        \RightLabel{$\texttt{Split}$}
        \UnaryInfC{$(p \rightarrow q_1) \vee (p \rightarrow q_2)$}
    \end{prooftree}
    is valid over $\mathfrak{C}$ relative to $\mathfrak{J} \cup \mathfrak{H}$. Without loss of generality, we can assume $\mathscr{D}_1$ to be closed and canonical, so we have
    \begin{prooftree}
        \AxiomC{$[p]_1$}
        \noLine
        \UnaryInfC{$\mathscr{D}_3$}
        \noLine
        \UnaryInfC{$q_1 \vee q_2$}
        \RightLabel{$1$}
        \UnaryInfC{$p \rightarrow q_1 \vee q_2$}
        \LeftLabel{$\mathscr{D}_2 = $}
        \RightLabel{$\texttt{Split}$}
        \UnaryInfC{$(p \rightarrow q_1) \vee (p \rightarrow q_2)$}
    \end{prooftree}
    We now distinguish two cases.
    \begin{enumerate}
        \item $\mathscr{D}_3$ is closed---namely, its conclusion $q_1 \vee q_2$ does not depend on assumption $p$, and $p$ is vacuously discharged by the step of $\rightarrow_I$ yielding the conclusion of $\mathscr{D}_1$. Since it must be $\langle \mathscr{D}_3, \mathfrak{J} \rangle$ valid on $\mathfrak{C}$, we have that $\mathscr{D}_3$ reduces modulo $\mathfrak{J}$ to an argument structure of the form
        \begin{prooftree}
            \AxiomC{$\mathscr{D}_4$}
            \noLine
            \UnaryInfC{$q_i$}
            \UnaryInfC{$q_1 \vee q_2$}
        \end{prooftree}
        with $\mathscr{D}_4 \in \texttt{DER}_\mathfrak{C}$, so that $\mathscr{D}_2 \leq_\mathfrak{J} \mathscr{D}_5$ where
        \begin{prooftree}
            \AxiomC{$\mathscr{D}_4$}
            \noLine
            \UnaryInfC{$q_i$}
            \UnaryInfC{$q_1 \vee q_2$}
            \UnaryInfC{$p \rightarrow q_1 \vee q_2$}
            \LeftLabel{$\mathscr{D}_5 = $}
            \RightLabel{$\texttt{Split}$}
            \UnaryInfC{$(p \rightarrow q_1) \vee (p \rightarrow q_2)$}
        \end{prooftree}
        We can hence apply $\phi_2$ and obtain
        \begin{prooftree}
            \AxiomC{$[p]_1$}
            \noLine
            \UnaryInfC{$\mathscr{D}^*_4$}
            \noLine
            \UnaryInfC{$q_i$}
            \RightLabel{$1$}
            \UnaryInfC{$p \rightarrow q_i$}
            \UnaryInfC{$(p \rightarrow q_1) \vee (p \rightarrow q_2)$}
        \end{prooftree}
        If axiom $p$ belongs to $\mathfrak{C}$, and is used (undischarged) in $\mathscr{D}_4$, axiom $p$ is turned into an assumption by $\phi_2$, and discharged in the second-to-last step of $\mathscr{D}^*_4$, where the latter is just $\mathscr{D}_4$ with the exception that horizontal bars on top of $p$ have been removed. By Corollary 3, $\langle \mathscr{D}^*_4, \mathfrak{J} \cup \mathfrak{H} \rangle$ is valid over $\mathfrak{C}$, hence $\langle \mathscr{D}_5, \mathfrak{J} \cup \mathfrak{H} \rangle$ is also so. The same obtains when axiom $p$ is not used (undischarged) in $\mathscr{D}_4$, as in this case $\mathscr{D}^*_4 = \mathscr{D}_4$. Since $\mathscr{D}_2 \leq_{\mathfrak{J} \cup \mathfrak{H}} \mathscr{D}_5$, we have that $\langle \mathscr{D}_2, \mathfrak{J} \cup \mathfrak{H}\rangle$ is also valid on $\mathfrak{C}$.
        \item $\mathscr{D}_3$ is open. Then we apply $\phi_1$ and obtain
        \begin{prooftree}
            \AxiomC{}
            \UnaryInfC{$p$}
            \noLine
            \UnaryInfC{$\mathscr{D}^*_3$}
            \noLine
            \UnaryInfC{$q_1 \vee q_2$}
            \UnaryInfC{$p \rightarrow q_1 \vee q_2$}
            \LeftLabel{$\mathscr{D}_5 = $}
            \RightLabel{$\texttt{Split}$}
            \UnaryInfC{$(p \rightarrow q_1) \vee (p \rightarrow q_2)$}
        \end{prooftree}
        Now, it must be $\langle \mathscr{D}_3, \mathfrak{J} \rangle$ valid on $\mathfrak{C}$. So, the closed instance $\mathscr{D}^*_3$ of $\mathscr{D}_3$ is valid relative to $\mathfrak{J}$ on the extension $\mathfrak{D}$ of $\mathfrak{C}$ obtained by adding axiom $p$ to $\mathfrak{C}$ (which is possibly equal to $\mathfrak{C}$). This instance reduces modulo $\mathfrak{J}$ to $\mathscr{D}_4$ as above. If axiom $p$ is not used (undischarged) in $\mathscr{D}_4$, then $\mathscr{D}_4 \in \texttt{DER}_\mathfrak{C}$, and we apply $\phi_2$ reasoning as in the case of point 1 above when it was assumed that axiom $p$ was not used (undischarged) in $\mathscr{D}_4$. If axiom $p$ is used (undischarged) in $\mathscr{D}_4$, then $\mathscr{D}_4 \in \texttt{DER}_\mathfrak{D}$. We apply $\phi_2$ and reason as in the case of point 1 above when it was assumed that axiom $p$ was used (undischarged) in $\mathscr{D}_4$.
    \end{enumerate}
    In all cases, hence, $\langle \mathscr{D}_2, \mathfrak{J} \cup \mathfrak{H} \rangle$ is valid on $\mathfrak{C}$, and the result now follows by arbitrariness of $\mathfrak{C}$ and of $\mathfrak{B}$.
\end{proof}

\begin{corollary}
    $p \rightarrow q \vee r \models (p \rightarrow q) \vee (p \rightarrow r)$
\end{corollary}

\begin{corollary}
    Conjecture 1 and Conjecture 2 fail.
\end{corollary}

\noindent Observe that Theorem 3 and Corollaries 4 and 5 are available as soon as one has atomic bases of level $0$, and independently of how $\bot$ is dealt with relative to the alternatives discussed in footnotes 2 and 4---i.e., with $\bot$ as a nullary connective whose validity on the base is tantamount to the validity of every atom on the base, rather than with $\bot$ as an atomic constant and bases which always contain all the atomic instances of \emph{ex falso}. Hence, the result applies to all variants of PTS with argument structures and reductions, where a constructivity and uniformity constraint is put on the latter.\footnote{The reductions $\phi^1$ and $\phi^2$ are deterministic only because, via Definition 5, atomic derivations are not allowed to contain or discharge assumptions. This feature was introduced for a better treatment of dischargement of rules of level $0$ by rules of level $2$, but is otherwise not essential. For, a proof similar to that of Theorem 3 can be found also when one allows atomic derivations to involve unbound assumptions. In this case we need two reductions, one which is defined precisely like $\phi^1$ above, while the other is defined on argument structures ending by $\texttt{Split}$ whose immediate sub-structure is canonical with immediate sub-structure closed and canonical, and which behaves like $\phi^2$ instead.}

\section{Conclusion: closure under substitutions and Pezlar function}

A first remark on Theorem 3 is that its proof relies essentially on a monotonic kind of validity. The proof goes through only if we know that a given argument structure from assumption $p$ to conclusion $q_i \ (i = 1, 2)$ or $q_1 \vee q_2$ yields, thanks to its monotonic validity relative to some $\mathfrak{J}$ on some $\mathfrak{B}$, an atomic derivation in $\mathfrak{B}$ plus $p$ as an axiom, when assumption $p$ is replaced by $p$ as an axiom. Thus, the proof \emph{does not} apply to non-monotonic P-tV.

A second, more important remark, is that the validity of the atomic Split rule is not closed under substitutions of atoms with formulas. I.e., we have not shown that, for every function $\star \colon \texttt{ATOM}_\mathscr{L} \to \texttt{FORM}_\mathscr{L}$,

\begin{prooftree}
    \AxiomC{$\star(p) \rightarrow \star(q) \vee \star(r)$}
    \UnaryInfC{$(\star(p) \rightarrow \star(q)) \vee (\star(p) \rightarrow \star(r))$}
\end{prooftree}
is logically valid. There are in fact counterexamples to this, e.g.,

\begin{prooftree}
    \AxiomC{$p \vee q \rightarrow p \vee q$}
    \UnaryInfC{$(p \vee q \rightarrow p) \vee (p \vee q \rightarrow q)$}
\end{prooftree}
can be shown \emph{not to be} logically valid in Prawitz's original miP-tV. Thus, as observed by \cite{piechaschroeder-heisterdecampossanz}, an alternative way of formulating incompleteness of $\texttt{IL}$ over $\models$ is by saying that $\models$ is not closed under substitutions.

As already remarked in Section 2, the incompleteness results proved by de Campos Sanz, Piecha and Schroeder-Heister for miB-eS, are much stronger than this, as they do not apply just to completeness for $\Vdash$ as per Definition 9, but to $\Vdash$ closed under substitutions, written $\Vdash^*$. If by $\star(A)$ we write the result of uniformly replacing atoms by formulas in $A$ through $\star$, and by $\star(\Gamma)$ we write the result of doing this for all the formulas of $\Gamma$, we have

\begin{center}
    $\Gamma \Vdash^* A$ iff $\forall \star (\star(\Gamma) \Vdash \star(A))$.
\end{center}
Then we put

\begin{center}
    $\texttt{IL}$ is $^*$-complete over $\Vdash^*$ iff $\Gamma \Vdash^* A \Longrightarrow \Gamma \vdash_{\texttt{IL}} A$
\end{center}
and we can prove that this fails in the miB-eS framework considered in this paper. Be that as it may, Prawitz's approach in \cite{prawitz1973}, and the conjecture that Prawitz formulates there, \emph{do not} officially require that logical validity of inference rules, and the relation of logical consequence $\models$ in terms of existence of logically valid arguments, hold under replacements of atoms with formulas. In other words, the kind of completeness at issue is the one provided by Definitions 26 and 27, which is why I took Theorem 3 to refute Conjectures 1 and 2.

Of course, one could require closure under substitutions also for Prawitz's approach, i.e., if we start with the miP-tV notion of consequence,

\begin{center}
    $\Gamma \models^* A$ iff $\forall \star (\star(\Gamma) \models \star(A))$
\end{center}
Observe that, because of how $\models$ is defined (Definition 24), $\Gamma \models^* A$ comes to mean that, for every $\star$, there is a logically valid $\langle \mathscr{D}, \mathfrak{J} \rangle$, with $\mathscr{D}$ from $\star(\Gamma)$ to $\star(A)$. This is fine, but can be strengthened. Rather than requiring the existence of a logically valid argument for each substitution-instance, we could require the existence of one single logically valid argument which remains so when the atoms in its formulas are replaced according to any substitution function. If we indicate by $\star(\mathscr{D})$ the result of such a replacement in $\mathscr{D}$, we would thus have

\begin{center}
    $\Gamma \models^* A$ iff $\exists \langle \mathscr{D}, \mathfrak{J} \rangle \forall \star (\langle \star(\mathscr{D}), \mathfrak{J} \rangle$ is logically valid, with $\star(\mathscr{D})$ from $\star(\Gamma)$ to $\star(A))$.
\end{center}
It is not clear which option is preferable, and the impression is that the choice depends mostly on personal (constructivist) tastes---compare also with the observation about the inversion in the order of the quantifiers in the definition of logical consequence discussed in Section 4.

However, if one decides to pursue the closure-under-substitution strategy, a source of inspiration might come from a recent result proved by Pezlar in \cite{pezlarselector} in a context more akin to Martin-L\"{o}f's type theory \cite{martin-loef}, via the formulas-as-types conception \cite{howard} and, hence, with rules enriched by constructors for given types. Pezlar focuses on a special form of the Split rule, i.e.,

\begin{prooftree}
    \AxiomC{$A \rightarrow B \vee C$}
    \RightLabel{$\texttt{Split}^*$}
    \UnaryInfC{$(A \rightarrow B) \vee (A \rightarrow C)$}
\end{prooftree}
where $A$ is a Harrop formula. $\texttt{Split}^*$ is in turn understood as a generalised elimination for $\vee$, i.e., as the rule

\begin{prooftree}
    \AxiomC{$[A]$}
    \noLine
    \UnaryInfC{$B \vee C$}
    \AxiomC{$[A \rightarrow B]$}
    \noLine
    \UnaryInfC{$D$}
    \AxiomC{$[A \rightarrow C]$}
    \noLine
    \UnaryInfC{$D$}
    \RightLabel{$\texttt{S}$}
    \TrinaryInfC{$D$}
\end{prooftree}
We now consider the following reductions:

\begin{footnotesize}
\begin{prooftree}
    \AxiomC{$\mathscr{D}$}
    \noLine
    \UnaryInfC{$A \rightarrow B \vee C$}
    \RightLabel{$\texttt{Split}^*$}
    \UnaryInfC{$(A \rightarrow B) \vee (A \rightarrow C)$}
    \AxiomC{$\stackrel{\texttt{Split}^*\texttt{-to-S}}{\Longrightarrow}$}
    \AxiomC{$\mathscr{D}$}
    \noLine
    \UnaryInfC{$A \rightarrow B \vee C$}
    \AxiomC{$[A]_1$}
    \BinaryInfC{$B \vee C$}
    \AxiomC{$[A \rightarrow B]_2$}
    \UnaryInfC{$(A \rightarrow B) \vee (A \rightarrow C)$}
    \AxiomC{$[A \rightarrow C]_3$}
    \UnaryInfC{$(A \rightarrow B) \vee (A \rightarrow C)$}
    \RightLabel{$\texttt{S}, 1, 2, 3$}
    \TrinaryInfC{$(A \rightarrow B) \vee (A \rightarrow C)$}
    \noLine
    \TrinaryInfC{}
\end{prooftree}
\end{footnotesize}

\begin{prooftree}
    \AxiomC{$[A]_1$}
    \noLine
    \UnaryInfC{$\mathscr{D}_1$}
    \noLine
    \UnaryInfC{$B_i$}
    \UnaryInfC{$B_1 \vee B_2$}
    \AxiomC{$[A \rightarrow B_1]_2$}
    \noLine
    \UnaryInfC{$\mathscr{D}_2$}
    \noLine
    \UnaryInfC{$D$}
    \AxiomC{$[A \rightarrow B_2]_3$}
    \noLine
    \UnaryInfC{$\mathscr{D}_3$}
    \noLine
    \UnaryInfC{$D$}
    \RightLabel{$\texttt{S}, 1, 2, 3$}
    \TrinaryInfC{$D$}
    \AxiomC{$\stackrel{\phi_{\texttt{S}}}{\Longrightarrow}$}
    \noLine
    \UnaryInfC{}
    \AxiomC{$[A]_1$}
    \noLine
    \UnaryInfC{$\mathscr{D}_1$}
    \noLine
    \UnaryInfC{$B_i$}
    \RightLabel{$1$}
    \UnaryInfC{$[A \rightarrow B_i]$}
    \noLine
    \UnaryInfC{$\mathscr{D}_{i + 1}$}
    \noLine
    \UnaryInfC{$D$}
    \noLine
    \TrinaryInfC{}
\end{prooftree}
Then, $\texttt{Split}^*$ is shown to be constructively and logically valid relative to the set $\{\texttt{Split}^*\texttt{-to-S}, \phi_\texttt{S}, \phi_\rightarrow\}$. The reason is essentially that, as proved by Smith in \cite{smith} using Kleene's Slash computability relation \cite{kleene1, kleene2}, an arbitrary construction for a Harrop formula will always yield a \emph{canonical} dependent construction for a corresponding dependent formula. It follows that the reduction $\phi_\texttt{S}$ is in this context in order, since it can be grounded on the following elimination rule in a kind of Martin-L\"{o}f style

\begin{prooftree}
    \AxiomC{$z : A$}
    \noLine
    \UnaryInfC{$c(z) : B \vee C$}
    \AxiomC{$[x : A \rightarrow B]$}
    \noLine
    \UnaryInfC{$d(x) : D$}
    \AxiomC{$[y : A \rightarrow C]$}
    \noLine
    \UnaryInfC{$e(y) : D$}
    \RightLabel{$\texttt{S}^\tau_1$}
    \TrinaryInfC{$\texttt{Pezlar}(z.c, x.d, y.e) : D$}
\end{prooftree}
where the Pezlar function is defined by the following equality rule, again in a kind of Martin-L\"{o}f style,

\begin{prooftree}
    \AxiomC{$[z : A]$}
    \noLine
    \UnaryInfC{$a(z) : B_i$}
    \AxiomC{$[x : A \rightarrow B_1]$}
    \noLine
    \UnaryInfC{$d_1(x) : D$}
    \AxiomC{$[y : A \rightarrow B_2]$}
    \noLine
    \UnaryInfC{$d_2(y) : D$}
    \RightLabel{$\texttt{S}^\tau_2$}
    \TrinaryInfC{$\texttt{Pezlar}(z.\texttt{in}_i(a(z)), x.d_1, y.d_2) = d_i(\lambda z.a(z)) : D$}    
\end{prooftree}
with $\texttt{in}_i$ injection function for forming a construction for $B_1 \vee B_2$ out of one for $B_i$. 

Pezlar's nice proof does not apply to Prawitz's framework, since here we cannot guarantee that any open argument structure from some Harrop formula valid on a given $\mathfrak{B}$ relative to some $\mathfrak{J}$ will reduce to an open \emph{canonical} argument structure from the same Harrop formula valid on $\mathfrak{B}$ relative to $\mathfrak{J}$. The reason is again that sets of reductions are a sub-class of the class of partial constructive functions---in fact, Smith's results are \emph{existence} theorems. However, observe also that the Pezlar function is constructive and uniform in a very strong sense, as it can be given by an explicit and purely schematic equation, and that, altough $\texttt{Split}^*$ is not closed under substitution either, it is of course much more general than $\texttt{Split}$. These topics can be discussed in future works.

To conclude, let me just remark that, while the closure-under-substitutions strategy is certainly available for miP-tV too, there might be independent reasons why a condition of closure under substitutions \emph{should not} be put on Prawitz's notion of validity. Closure under substitutions is a crucial property with respect to a \emph{primitive} notion of proof-theoretic consequence, as the one of miB-eS is. However, when the notion is \emph{not} primitive, but defined in terms of existence of a proof from a given set of assumptions to a given conclusion, as happens in miP-tV, one is not just focusing on \emph{what} one can show to be logically valid, but on \emph{how} such a logical validity can be established. This is a semantics of proofs, as opposed to a semantics of formulas.

When investigating intuitionistic completeness, hence, we may well ask whether everything that is proof-theoretically valid is also derivable in intuitionistic logic. But since our semantics is a semantics of \emph{proofs}, the question of paramount importance is rather whether any \emph{piece} of proof-theoretically valid reasoning can be safely derived in the chosen logic. And this piece of reasoning needs not be closed under substitutions, of course.

\paragraph{Acknowledgements} I am indebted to Hermógenes Oliveira, Ivo Pezlar, Thomas Piecha, Luca Tranchini, and especially Dag Prawitz, for discussions and insights on the content of this paper. This work has been supported by the grant PI 1965/1-1 for the DFG project \emph{Revolutions and paradigms in logic. The case of proof-theoretic semantics}.

\paragraph{Conflict of interests} The author declares that there is no conflict of interests.

\bibliographystyle{abbrv}
\bibliography{bibliography}

\end{document}